\documentstyle[12pt]{article}
\font\sym=msbm10 scaled \magstep1
\newcommand{\IZ}{\hbox{\sym \char '132}}

\newcommand{\IP}{\hbox{\sym \char '120}}

\newcommand{\IH}{\hbox{\sym \char '110}}
\newcommand{\IC}{\hbox{\sym \char '103}}

\newcommand{\zz}{\IZ}

\newcommand{\ke}{{\cal E}}
\newcommand{\kf}{{\cal F}}
\newcommand{\kg}{{\cal G}}
\newcommand{\kh}{{\cal H}}
\newcommand{\ki}{{\cal I}}
\newcommand{\kk}{{\cal K}}
\newcommand{\kl}{{\cal L}}
\newcommand{\km}{{\cal M}}
\newcommand{\kn}{{\cal N}}
\newcommand{\ko}{{\cal O}}

\newcommand{\kq}{{\cal Q}}
\newcommand{\ks}{{\cal S}}

\newcommand{\ku}{{\cal U}}

\newcommand{\bP}{\IP}

\newcommand{\kea}{{\ke_{\alpha}}}

\newcommand{\atilde}{{\tilde a}}
\newcommand{\qtilde}{{\tilde q}}
\newcommand{\ketilde}{{\tilde\ke}}

\newcommand{\ph}{\varphi}
\newcommand{\eps}{\varepsilon}
\newcommand{\taus}{\tau^{ss}}
\newcommand{\vb}{vector bundle~}
\newcommand{\vbs}{vector bundles~}
\newcommand{\ms}{moduli space~}
\newcommand{\mss}{moduli spaces~}
\newcommand{\st}{stability~}
\newcommand{\lra}{\longrightarrow}
\newcommand{\isom}{\cong}
\newcommand{\dual}{^{\rm v}}
\newcommand{\ddual}{^{\rm vv}}

\newcommand{\Sl}{{\rm SL}}
\newcommand{\Gl}{{\rm GL}}
\newcommand{\GL}{{\rm GL}}
\newcommand{\PGL}{{\rm PGL}}
\newcommand{\Gm}{{\rm\bf G}_m}
\newcommand{\rk}{{\rm rk}}
\newcommand{\Ker}{{\rm Ker }}
\newcommand{\coker}{{\rm coker }}
\newcommand{\Ta}{{([a],[T])}}
\newcommand{\Tors}{{\rm T}}
\newcommand{\Quot}{{\rm Quot }}
\newcommand{\Pic}{{\rm Pic }}
\newcommand{\Ima}{{\rm Im}}
\newcommand{\Hom}{{\rm Hom }}
\newcommand{\id}{{\rm id }}

\newcommand{\Supp}{{\rm Supp}}
\newcommand{\Spec}{{\rm Spec}}
\newtheorem{theorem}{Theorem}[section]
\newtheorem{proposition}[theorem]{Proposition}

\newtheorem{lemma}[theorem]{Lemma}
\newtheorem{corollary}[theorem]{Corollary}
\newtheorem{remark}[theorem]{Remark}
\newtheorem{definition}[theorem]{Definition}

\newcommand{\prf}{{\it Proof:} }

\newcommand{\qed}{\hspace*{\fill}\hbox{$\Box$}}
\newcommand{\ses}[3]{0\longrightarrow#1\longrightarrow#2\longrightarrow#3
   \longrightarrow0}
\newcommand{\sesq}[5]{0\lra#1\stackrel{#2}{\lra}#3\stackrel{#4}{\lra}#5\lra0}
\newcommand{\rpfeil}[2]{\stackrel{\raisebox{-1mm}{$\scriptstyle#2$}}{{\hbox
to #1mm{\rightarrowfill}}}}

\newcommand{\spfeil}[2]{{\scriptstyle{#2}\downarrow\phantom{{\scriptstyle#2}}}}
\newcommand{\epimorph}{\longrightarrow\!\!\!\!\!\!\!\!\longrightarrow}
\newcommand{\btimes}{\hbox{$\times\!\!\>\!\!\!\!\lower1.25pt\hbox{$\Box$}\,$}}

\renewcommand{\.}{{\cdot}}

\newcommand{\propsubset}{\, \mbox{\raisebox{1mm}{\makebox[0mm][c]{\,$\subset$}}
\raisebox{-1.5mm}{\makebox[0mm][c]{$\neq$\,\,\,}}}\,}

\begin{document}
\title{Stable pairs on curves and surfaces}
\author{D. Huybrechts\thanks{Max-Planck-Institut fur Mathematik,
Gottfried-Claren-Str. 26, 5300 Bonn 3, Germany}
\and
M. Lehn\thanks{Math. Institut der Univ., R\"amistr. 74, 8001 Z\"urich,
Switzerland}}
\date{November 1992}
\maketitle

{\large \bf Introduction}

During the last years there has been growing interest
in vector bundles with additional
structures, e.g. parabolic and level structures.
This paper results from an attempt to construct quasi-projective moduli
spaces for
framed bundles, i.e. bundles together with an isomorphism to a fixed bundle
on a divisor
as introduced in \cite{d2}, \cite{l4} and \cite{l2}. More generally one can
ask for bundles with a homomorphism
to a fixed sheaf $\ke_0$. We use techniques of geometric invariant
theory to construct projective moduli spaces. This leads to natural
stability conditions. In contrast to the pure
bundle case an extra parameter appears in the definition of stability.

A pair $(\ke,\alpha)$ consisting of a coherent sheaf $\ke$
 on a smooth, projective variety and
a homomorphism $\alpha$ from $\ke$ to $\ke_0$ is called {\it stable}
with respect to a polynomial $\delta$
if and only if the following conditions are satisfied.

{\it i)} $\chi_\kg<{(\rk\kg/\rk\ke)}\chi_\ke-{(\rk \kg/\rk\ke)}\delta$ for
all subsheaves
$\kg\subset\Ker\alpha$.

{\it ii)} $\chi_\kg< {(\rk\kg/\rk \ke)}\chi_\ke+\delta{(\rk\ke-\rk\kg)/\rk\ke}
$ for all
subsheaves $\kg\,\propsubset\ke$.

Here $\chi$ denotes the Hilbert polynomial and the inequalities must
hold for large arguments.
In \S 1 we prove

{\bf Theorem:} {\it For a smooth, projective variety $X$ of dimension one or
two there is a fine quasi-projective moduli
space of stable pairs $(\ke,\alpha:\ke\to\ke_0)$ with respect to $\delta$.}

Moreover, we will prove that this space can be naturally compactified
(For a precise statement see \ref{mainthm}).

In particular, this theorem proves the quasi-projectivity of many of the \mss
of framed
bundles, which in \cite{l4} were constructed only as algebraic spaces
(\ref{framedareproj}).
In \S2 we study two
special cases for $\ke_0$, where $\ke_0$ is the structure sheaf $\ko_X$ or a
\vb on an effective divisor.

The case $\ke_0\cong\ko_X$ leads to the definition of Higgs pairs, i.e.
solution
of the vortex equation as considered in \cite{Br}, \cite{Be}, \cite{GP},
\cite{Th}. A Higgs pair is a vector bundle $\ke$ together with a global section
$\varphi$ satisfying certain stability conditions. The corresponding moduli
spaces of rank two vector bundles on a curve were constructed by M. Thaddeus
and A. Bertram.
Dualizing the situation one gets a vector bundle $\ke\dual$ together with a
homomorphism
$\alpha=\varphi\dual:\ke\dual\lra\ko_X$. The stability conditions for Higgs
pairs
 translate into {\it i)} and {\it ii)} above. This dual point of view
allows us to compactify the moduli space in the surface case, too, by adding
pairs with
torsionfree sheaves. Instead of one \ms M. Thaddeus consideres the whole
series of moduli spaces, which result from changing the stability parameter
in order to
'approximate' the usual \ms of semistable bundles. We generalize this method
for bundles on
a surface and describe the 'limit' of this series. As a generalization of
Bogomolov's result we
prove a theorem about the restriction of stable pairs
to curves of high degree (\ref{restofpairs}).

The case of $\ke_0$ being a \vb on a divisor leads to the concept of bundles
with
level structure (\cite{Se}) and to the concept of framed bundles (\cite{l4})
in
dimension one and two, resp.

\section{Moduli spaces of stable pairs}

Throughout this paper we fix the following notations: $X$ is
an irreducible,
nonsingular, projective variety of dimension $e$ over an
algebraically closed
field $k$ of characteristic zero, embedded by a very ample line bundle
$\ko_X(1)$.
The canonical
line bundle is denoted by $\kk_X$. If $\ke$ is a coherent $\ko_X$-module,
then
$\chi_{\ke}(n):=\chi(\ke\otimes\ko_X(n))$ denotes
its Hilbert polynomial,
$\Tors(\ke)$ its torsion submodule and $\det\,\ke$ its
determinant line bundle.
The degree of $\ke$, $\deg\,\ke$, is the integral
number $c_1(\det\,\ke).H^{e-1}$,
where $H\in|\ko_X(1)|$ is a hyperplane
section.

$\chi$ will always be a polynomial
with rational coefficients which
has the form
$$\chi(z)=\deg\,X\cdot
r\cdot\frac{z^e}{e!}+(d-\frac{\deg\,\kk_X}{2}\cdot
r)\cdot z^{e-1}+\hbox{\rm
Terms of lower order in $z$}.$$
If $\chi=\chi_{\ke}$, then $r=\rk\ke$ and $d=\deg\ke$.
Finally, let $\ke_0$ be
a fixed coherent $\ko_X$-module. By a \em pair \em we will
always mean a pair
$(\ke,\alpha)$ consisting of a coherent $\ko_X$-module $\ke$ with
Hilbert polynomial $\chi_{\ke}=\chi$ and a nontrivial homomorphism
$\alpha:\ke\rightarrow\ke_0$. We write $\kea$ for $\Ker\,\alpha$.

In the next section we define the notion of semistability for such pairs with
respect to an additional parameter $\delta$. To simplify the
notations and to be able to treat stability and semistability simultaneously,
we employ the following short-hand: Whenever in a statement the word \em
(semi)stable \em occurs together with a relation symbol in brackets, say
$(\leq)$, the latter should be read as $\leq$ in the semistable case and as
$<$ in the stable case. An inequality $p\,(\leq)\,p'$ between polynomials
means, that $p(n)\,(\leq)\,p'(n)$ for large integers $n$. If $p$ is a
polynomial then $\Delta p(n):=p(n)-p(n-1)$ is the difference polynomial.

We proceed as follows: In section 1.1 we define semistability for pairs and
formulate the moduli problem. In section 1.2 boundedness
results for semistable pairs on curves and surfaces are obtained. Moreover, a
close
relation between semistability and sectional semistability is established.
The notion of sectional stability naturally appears by way of constructing
moduli spaces
for pairs. This is done in section 1.3 leading to the existence theorem
\ref{mainthm}.
Section 1.4 is devoted to an invariant theoretical analysis of the
construction in 1.3 and the proof of the main technical proposition
\ref{prop}.\\
The reader who is familiar with the papers of Gieseker and Maruyama
(\cite{g2}, \cite{Ma})
will notice that many of our arguments are generalizations of their techniques.


\subsection{Stable pairs and the moduli problem}\label{Kapitel11}
Let $\delta$ be a polynomial with rational coefficients such that
$\delta>0$, i.\ e.\ $\delta(n)>0$ for all $n\gg 0$. We write
$\delta(z)=\sum_{\nu}\delta_{e-\nu}z^{\nu}$.

\begin{definition}\label{stabdef} A pair
$(\ke,\alpha)$ is called \em(semi)stable \em (with respect to $\delta$), if
the following two conditions are satisfied:
\begin{itemize}\item[(1)]
$\rk\ke\.\chi_{\kg}\,(\leq)\,\rk\kg\.(\chi_{\ke}-\delta)$ for all
nontrivial submodules $\kg\subseteq\kea$.
\item[(2)]
$\rk\ke\.\chi_{\kg}\,(\leq)\,\rk\kg\.(\chi_{\ke}-\delta)+\rk\ke\.\delta$ for
all nontrivial submodules $\kg\,(\subseteq)\,\ke$.
\end{itemize}\end{definition}

If no confusion can arise, we omit $\delta$ in the notations.
Note that a stable pair a fortiori is semistable.

\begin{lemma}\label{firstobs} Suppose $(\ke,\alpha)$ is a semistable pair,
then:
\begin{itemize}
\item[i)] $\kea$ is torsion free. $h^0(\kg)\leq h^0(\Tors(\ke_0))$ for
all submodules $\kg\subseteq\Tors(\ke)$.
\item[ii)] Unless $\alpha$ is injective, $\delta$ is a polynomial of degree
smaller than $d$.
\end{itemize}\end{lemma}

\prf ad i): If $\kg\subset\kea$ is torsion, then $\rk\kg=0$. Condition
(1) then shows $\chi_{\kg}=0$, hence $\kg=0$. Thus $\alpha$ embeds the torsion
of $\ke$ into the torsion of $\ke_0$. This gives the second assertion. ad ii):
Assume $\kea$ is nontrivial. By i) $\kea$ is torsion free of
positive rank, and condition (1) implies
$\delta/\rk\ke\,\leq\,(\chi_{\ke}/\rk\ke-\chi_{\kea}/\rk\kea)$. The two
fractions in the brackets are polynomials with the same leading coefficients.
This shows $\deg\delta<e$.\qed

Thus if $\deg\delta\geq e$, then $\alpha$ must needs be an injective
homomorphism, and isomorphism classes of semistable pairs
correspond to submodules of $\ke_0$ with fixed Hilbert polynomial.
Note that condition (2) of the definition above is automatically satisfied.
So in this case all pairs are in fact stable and parametrized by
the projective quotient
scheme $\Quot^{\chi_{\ke_0}-\chi_{\ke}}_{X/\ke_0}$. For that reason we assume
henceforth that $\delta$ has the form
$$\delta(z)=\delta_1z^{e-1}+\delta_2z^{e-2}+\cdots+\delta_e.$$

\begin{definition}\label{mustab} A pair $(\ke,\alpha)$ is called
$\mu$-(semi)stable
(with respect to $\delta_1$), if the following two conditions are satisfied:
\begin{itemize}
\item[(1)] $\rk\ke\.\deg\kg\,(\leq)\,\rk\kg\.(\deg\ke-\delta_1)$ for all
nontrivial submodules $\kg\subseteq\kea$.
\item[(2)]
$\rk\ke\.\deg\kg\,(\leq)\,\rk\kg\.(\deg\ke-\delta_1)+\rk\ke\.\delta_1$ for all
nontrivial submodules $\kg\subseteq\ke$ with $\rk\kg<\rk\ke$.
\end{itemize}\end{definition}

As in the theory of stable sheaves there are immediate implications for pairs
$(\ke,\alpha)$:
\begin{center}
$\mu$-stable $\Rightarrow$ stable $\Rightarrow$ semistable $\Rightarrow$
$\mu$-semistable
\end{center}

A family of pairs parametrized by a noetherian scheme $T$ consists of a
coherent $\ko_{T\times X}$-module $\ke$, which is flat over $T$, and a
homomorphism $\alpha:\ke\rightarrow p_X^*\ke_0$. If $t$ is a point of $T$, let
$X_t$ denote the fibre $X\times\Spec\,k(t)$, $\ke_t$ and $\alpha_t$ the
restrictions of $\ke$ and $\alpha$ to $X_t$. A homomorphism of pairs
$\Phi:(\ke,\alpha)\rightarrow(\ke',\alpha')$ is a homomorphism
$\Phi:\ke\rightarrow\ke'$ which commutes with $\alpha$ and $\alpha'$, i.\ e.\
$\alpha'\circ\Phi=\alpha$. The correspondence
$$T\mapsto\{\hbox{Isomorphism classes of families of (semi)stable pairs
parametrized by
$T$}\}$$ defines a setvalued contravariant functor
$\underline\km_\delta^{(s)s}(\chi,\ke_0)$
on the category of noetherian
$k$-schemes of finite type. We will prove that for $\dim X\leq 2$
there is a fine moduli space for $\underline
\km_\delta^s
(\chi,\ke_0)$. It is compactified by equivalence classes of semistable pairs
(\ref{mainthm}).


\subsection{Boundedness and sectional stability}\label{Kapitel12}
In section \ref{Kapitel13} we will construct moduli spaces of stable pairs by
means of geometric invariant theory. The stability property needed in this
construction differs slightly from the one given in \ref{stabdef} in refering
to the number of global sections rather than to the Euler characteristic of a
submodule of $\ke$. In this section we compare the different notions and prove
that semistable pairs form bounded families, if the variety $X$ is a curve or a
surface.

\begin{definition}\label{secstadef} Let $\bar\delta$ be a positive rational
number. A pair $(\ke,\alpha)$ is called \em sectional (semi)\-stable \em (with
respect to $\bar\delta$), if $\kea$ is torsionfree and there is a subspace
$V\subseteq H^0(\ke)$ of dimension $\chi(\ke)$ such that the following
conditions are satisfied:
\begin{itemize}
\item[(1)] $\rk\ke\.\dim(H^0(\kg)\cap
V)\,(\leq)\,\rk\kg\.(\chi(\ke)-\bar\delta)$
for all nontrivial submodules $\kg\subseteq\kea$.
\item[(2)]
$\rk\ke\cdot\dim(H^0(\kg)\cap
V)\,(\leq)\,\rk\kg\cdot(\chi(\ke)-\bar\delta)+\rk\ke\cdot\bar\delta$
for all nontrivial submodules $\kg\,(\subseteq)\,\ke$.
\end{itemize}
\end{definition}

We begin with the case of a curve. In this case $\delta$ is a
rational number, and the Hilbert polynomial of any $\ko_X$-module $\kg$
depends on $\rk\kg$ and $\deg\kg$ only. Moreover, the polynomials occuring
in the inequalities of definition \ref{stabdef} are linear and have the same
leading coefficients. Therefore the Hilbert polynomials $\chi_{\kg}$ can
throughout be replaced by the Euler characteristics $\chi(\kg)$ without
changing the essence of the definition.

\begin{theorem}\label{boundcurve} Let $X$ be a smooth curve of genus $g$.
Assume that $d>r\.(2g-1)+\delta$.
\begin{itemize}\item[i)] If $(\ke,\alpha)$ is semistable or
sectional semistable, then $\ke$ is globally generated and $h^1(\ke)=0$.
\item[ii)] $(\ke,\alpha)$ is a (semi)stable pair if and only if it is
sectional (semi)stable.\end{itemize}
\end{theorem}

\prf ad i): On a smooth curve $X$ there is a split short exact sequence
$$\ses{\Tors(\ke)}{\ke}{\bar\ke}$$
with locally free $\bar\ke$ for any coherent $\ko_X$-module $\ke$. Now
$H^1(\ke)=H^1(\bar\ke)$, and $\ke$ is globally generated if and only if
$\bar\ke$ is globally generated. A glance at the short exact sequence
$$\ses{\bar\ke(-x)}{\bar\ke}{\bar\ke\otimes\ko_x}$$
for some closed point $x\in X$  shows that the vanishing of $H^1(\bar\ke(-x))$
for all $x\in X$ is a sufficient criterion for both $H^1(\ke)=0$ and the
global generation of
$\ke$. If $H^1(\bar\ke(-x))\neq0$, then there is a nontrivial homomorphism
$\ph:\bar\ke\rightarrow\kk_X(x)$. Let $\kg:=\Tors(\ke)+\Ker\ph$, so that there
is a short exact sequence
$$\ses{\kg}{\ke}{\kk_X(x-C)}$$
with some effective divisor $C$ on $X$. From this sequence we get
$$\chi(\kg)\geq\chi(\ke)-\chi(\kk_X(x))\quad{\rm
and}\quad h^0(\kg)\geq h^0(\ke)-h^0(\kk_X(x)).$$
On the other hand,
$$\chi(\kg)\leq\frac{\rk\ke-1}{\rk\ke}\chi(\ke)+\frac{\delta}{\rk\ke},$$
if $(\ke,\alpha)$ is semistable, and
$$\dim(V\cap H^0(\kg))\leq\frac{\rk\ke-1}{\rk\ke}\chi(\ke)
+\frac{\delta}{\rk\ke}$$
for some vector space $V\subseteq H^0(\ke)$ of dimension $\chi(\ke)$, if
$(\ke,\alpha)$ is sectional semistable. In the first case we get
$\chi(\ke)\leq\rk\ke\.\chi(\kk_X(x))+\delta$. And in the second
case one has
\begin{eqnarray*}h^0(\ke)-h^0(\kk_X(x))\leq h^0(\kg)&\leq&
\dim(H^0(\kg)\cap V)+(h^0(\ke)-\dim V)\\
&\leq&\frac{\rk\ke-1}{\rk\ke}\.\chi(\ke)+\frac{\delta}{\rk\ke}+(h^0(\ke)-
\chi(\ke))
\end{eqnarray*}
So in any case we end up with
$\deg\ke\leq\rk\ke\.(2g-1)+\delta$ contradicting the assumption of the
theorem.

ad ii): By part i) we have $\chi(\ke)=h^0(\ke)$, $V=H^0(\ke)$ and, of course,
$\chi(\kg)\leq h^0(\kg)$ for any submodule $\kg\subseteq\ke$. Hence
sectional (semi)stability
implies (semi)stability at once. Conversely, assume that $(\ke,\alpha)$ is a
(semi)stable pair. If for a submodule $\kg$ we have $h^1(\kg)=0$, then
$h^0(\kg)=\chi(\kg)$ and there is nothing to show. (This applies in particular
when
$\rk\kg=0$). Hence assume $h^1(\kg)\neq0$. As above this leads to a short
exact sequence
$$\ses{\kg'}{\kg}{\kk_X(-C)}$$
with $\rk\kg'=\rk\kg-1$ and some effective divisor $C$ on $X$, so that
$h^0(\kg')\geq h^0(\kg)-g$. By induction we may assume that
$$h^0(\kg')\,(\leq)\,\frac{\rk\kg-1}{\rk\ke}(h^0(\ke)-\delta)+\eps\.\delta$$
with $\eps=0$ if $\kg\subseteq\kea$ and $\eps=1$ if $\kg\,(\subseteq)\,\ke$.
Combining these inequalities we get
$$h^0(\kg)\,(\leq)\,\frac{\rk\kg}{\rk\ke}(h^0(\ke)-\delta)+
\eps\.\delta+(g-\frac{h^0(\ke)}{\rk\ke}+\frac{\delta}{\rk\ke}).$$
Since $h^0(\ke)=\chi(\ke)=\deg\ke+(1-g)\rk\ke>g\.\rk\ke+\delta$, we are
done.\qed

\begin{corollary} Suppose $X$ is a curve. The set of isomorphism classes of
$\ko_X$-modules occuring in semistable pairs is bounded.\qed
\end{corollary}

Before we pass on to surfaces recall the following criterion due to Kleiman
which we will use several times:

\begin{theorem}[Boundedness criterion of Kleiman]\label{boundcrit}
Suppose $\chi$ is a polynomial\\ and $K$ an integer. If $T$ is a set of
$\ko_X$-modules $\kf$ such that $\chi_{\kf}=\chi$ and
$$h^0(X,\kf|_{H_1\cap\ldots\cap H_i})\leq K\qquad\forall\,i=0,\ldots,e,$$
for a $\kf$-regular sequence of hyperplane sections $H_1,\ldots,H_e$, then
$T$ is bounded.\end{theorem}

\prf \cite[Thm 1.13]{k1}\qed

We introduce the following notation: For integers $\rho$ and $\eps$ let
$P(\rho,\eps)$ be the polynomial
$$P(\rho,\eps,z):=\frac{\rho}{r}(\chi(z)-\delta(z))+\eps\.\delta(z).$$
If $\kg\subseteq\ke$ is a submodule, let $\eps(\kg)=0$ or $1$ depending on
wether $\kg\subseteq\kea$ or not. Then the stability conditions can be
conveniently reformulated:
\begin{itemize}
\item[-] $(\ke,\alpha)$ is (semi)stable if and only if
$\chi_{\kg}\,(\leq)\,P(\rk\kg,\eps(\kg))$
for all nontrivial submodules $\kg\,(\subseteq)\,\ke$.
\item[-] $(\ke,\alpha)$ is $\mu$-(semi)stable if and only if
$\Delta\chi_{\kg}\,(\leq)\,\Delta P(\rk\kg,\eps(\kg))$ for all nontrivial
submodules $\kg\subseteq\ke$ with $\rk\kg<\rk\ke$.
\item[-] $(\ke,\alpha)$ is sectional (semi)stable if and only if
$\Tors(\kea)=0$ and there is a
subspace $V\subseteq H^0(\ke)$ of dimension $\chi(\ke)$ such that $\dim(V\cap
H^0(\kg))\,(\leq)\,P(\rk\kg,\eps(\kg),0)$ for all nontrivial submodules
$\kg\,(\subseteq)\,\ke$.
\end{itemize}

\begin{lemma}\label{n0} Suppose $X$ is a surface. There is an integer $n_0<0$,
depending on $X$, $\ko_X(1)$ and $P$ only, such that
$\Delta\chi_{\ko_X(-n_0)}>\Delta P(1,\eps)$ for $\eps=0,1$.
\end{lemma}

\prf As polynomials in $\nu$ the expressions $\Delta\chi_{\ko_X}(\nu-n)$ and
$\Delta P(1,\eps,\nu)$ are both linear and have the same positive leading
coefficient. Hence for very negative numbers $n$ one has
$\Delta\chi_{\ko_X}(\nu-n)>\Delta P(1,\eps,\nu)$.\qed

The following technical lemma is an adaptation of \cite[Lemma 1.2]{g2}.
Unfortunately, we cannot apply Gieseker's lemma directly because it treats
torsion free modules only, even though the necessary modifications are minor.

\begin{lemma}\label{techlem}
Suppose $X$ is a surface. Let $Q$ be a positive integer. Then there are
integers $N$ and $M$, depending on $X,\ko_X(1),P$ and $Q$,
such that if $\eps\in\{0,1\}$ and if $\kf$ is an $\ko_X$-module of rank
$r'\leq r$ with the properties
$h^0(\Tors(\kf))\leq Q$ and $\Delta\chi_{\kg}\leq\Delta P(\rk\kg,\eps)$ for
all nontrivial submodules $\kg\subseteq\kf$, then either
\begin{itemize}\item[] $h^0(\kf(n))<P(r',\eps,n)$ for all $n\geq N$,
\end{itemize}
or the following assertions hold:
\begin{itemize}
\item[(1)] $\Delta\chi_{\kf}=\Delta P(r',\eps)$,
\item[(2)] $h^2(\kf(n))=0$ for all $n\geq N$,
\item[(3)] $h^0(\kf(n_0)|_H)\leq M$ for some $\kf$-regular hyperplane section
$H$,
\item[(4)] if $h^1(\kf(n_0))\leq Q$, then $h^1(\kf(n))=0$ for all $n\geq N$.
\end{itemize}\end{lemma}

\prf Let $\kf$ be an $\ko_X$-module satisfying the assumptions of the lemma.
For every integer $n$ let $\kh_n'$ denote the image of the evaluation map
$H^0(\kf(n))\otimes\ko_X\rightarrow\kf(n)$ and $\ks_n'$ the quotient
$\kf(n)/\kh'_n$. Let $\kh_n$ be the kernel of the epimorphism
$$\kf(n)\epimorph (\ks_n'/\Tors(\ks_n'))=:\ks_n.$$
Then $\kh_n$ is characterized
by the following properties: $H^0(\kh_n)=H^0(\kf(n))$, $\kf(n)/\kh_n$ is
torsion free and $\kh_n$ is minimal with these two properties. Obviously
$\kh'_{n-1}(1)\subseteq\kh'_n$ and therefore also $\kh_{n-1}(1)\subseteq\kh_n$.
Moreover, being a submodule of the torsion free module $\kf(n-1)/\kh_{n-1}$
the quotient $\kh_n(-1)/\kh_{n-1}$ is itself torsion free. In particular
either
$\kh_{n-1}=\kh_n(-1)$ or $\rk\kh_{n-1}<\rk\kh_n$. Let $n_1<\ldots<n_k$ be the
indices with $\rk\kh_{n_i-1}<\rk\kh_{n_i}$. (If $\kf$ is torsion, then
$\kh_n=\kf(n)$ for all $n$. Let $k=0$ in this case). By Serre's Theorem
$\kh_{n_k}=\kf(n_k)$ and $k\leq r'$.

Let $s\in H^0\kf(n)$ be a nonzero section. Then either $s$ is a torsion
element or induces an injection $\ko_X(-n)\rightarrow\kf$. In the latter case
one has
$\Delta\chi_{\ko_X(-n)}\leq\Delta P(1,\eps)$. This is impossible for $n\leq
n_0$. It follows that
$$h^0(\kf(n_0))=h^0(\Tors(\kf)(n_0))\leq h^0(\Tors(\kf))\leq Q$$
and that $\kh_{n_0}=\Tors(\kf)(n_0)$. In particular $n_0<n_1$ if $r'>0$.

A generic hyperplane section
$H\in|\ko_X(1)|$ has the following properties:\begin{itemize}
\item[a)] $H$ is a smooth curve (of genus $g=1+\deg\kk_X/2)$.
\item[b)] $H$ is $\kh_n$-regular for all integers $n$.
\item[c)] $\kh_n|_H$ is globally generated at the generic point of $H$ for all
integers
$n$.
\end{itemize}
(a) is just Bertini's Theorem. For (b) it is enough to consider the sheaves
$\kh_{n_i}$, $i=0,\ldots,k$. $H$ must not contain any of the finitely many
associated points of the modules $\kh_{n_i}$ in the scheme $X$. But this is an
open condition. $\kh_n$ is globally generated outside the
support of $\Tors(\ks_n)$, so for (c) it is sufficient that in
addition $H$ should not contain any of the associated points of the
$\Tors(\ks_{n_i})$.
Hence for a generic hyperplane section $H$ there are short exact sequences
$$\ses{\kh_n(-1)}{\kh_n}{\kh_n|_H},$$
$$\ses{\ko_H^{r_n}}{\kh_n|_H}{Q_n},$$
where $r_n=\rk\kh_n$ and $Q_n$ is an $\ko_H$-torsion module. From the second
sequence one deduces estimates
$$h^1(\kh_n|_H)\leq r_n\.g\qquad{\rm and}\qquad h^1(\kh_n(\ell)|_H)=0,$$
if $\deg(K_H-\ell H)<0$, i.\ e.\ if $\ell>(2g-2)/H^2$. In particular we get
for all integers $n$ with $n_i+(2g-2)/H^2<n<n_{i+1}$:
$$h^1(\kh_n|_H)=h^1(\kh_{n_i}(n-n_i)|_H)=0.$$
This leads to the inequalities
$$\begin{array}{rcl}
h^0(\kf(n))-h^0(\kf(n-1))&=&h^0(\kh_n)-h^0(\kh_n(-1))\\
&\leq&h^0(\kh_n|_H)=\chi(\kh_n|_H)+h^1(\kh_n|_H)
\end{array}$$
and, summing up,
$$ h^0(\kf(n))-h^0(\kf(n_0)\leq
\sum_{\nu=n_0+1}^n\chi(\kh_{\nu}|_H)+\sum_{\rho=1}^{r_n}
\rho g((2g-2)/H^2+1).$$
Let $K:=Q+{r'+1\choose2}g((2g-2)/H^2+1)$. Then
$$h^0(\kf(n))\leq K+\sum_{\nu=n_0+1}^n\chi(\kh_{\nu}|_H)$$
for all integers $n\geq n_0$.
Suppose $n_0\leq\nu<n_k$. Then $r_{\nu}<r'$. Since $\kh_{\nu}(-\nu)$ is a
submodule of $\kf$,
$$
\chi(\kh_{\nu}|_H)=\Delta\chi_{\kh_{\nu}(-\nu)}(\nu)
\leq\Delta P(r_{\nu},\eps,\nu)$$
Now
\begin{eqnarray*}
\Delta P(r_{\nu},\eps,\nu)-\Delta P(r',\eps,\nu)&=&(r_{\nu}-r')\.(\deg
X\.\nu+d/r+(1-g)-\delta_1)\\
&\leq&-(\deg X\.\nu+C),
\end{eqnarray*}
where $C$ is a constant depending on $r$, $d$, $n_0$, $\deg X$ and $g$.
For $\nu\geq n_k$ one has $\kh_{\nu}=\kf(\nu)$ so that
$\chi(\kh_{\nu}|_H)=\Delta\chi_{\kf}(\nu)$. Let $m(n)=\min\{n,n_k-1\}$. Then
the following inequality holds for all $n\geq n_0$:
\begin{eqnarray*}
h^0(\kf(n))-\sum_{\nu=n_0+1}^{n}\Delta P(r',\eps,\nu) &\leq&
K-\sum_{\nu=n_0+1}^{m(n)}\{\deg X\.\nu+C\}\\
&&-\sum_{\nu=m(n)+1}^{n}\{\Delta P(r',\eps,\nu)-\Delta\chi_{\kf}(\nu)\}.
\end{eqnarray*}
Note that the summands of the second sum of the right hand side are all equal
to some nonnegative constant $C'$, (and that by convention the sum is 0 if
$n<n_k$). Let $f$ be the polynomial
$$f(z):=\deg X({z+1\choose 2}-{n_0+1\choose 2})+C\.(z-n_0)-K-P(r',\eps,n_0).$$
Then for $n\geq n_0$:
$$h^0(\kf(n))-P(r',\eps,n)\leq-f(m(n))-C'\.(n-m(n))$$
There is an integer $N_1>n_0$ such that $f(\nu)>0$ for all
$\nu\geq N_1$. Assume $N>N_1$. If $n_k-1\geq N_1$ then for all $n\geq N$ one
has $m(n)\geq N_1$, hence $f(m(n))>0$ and $h^0(\kf(n))<P(r',\eps,n)$.
Hence we can restrict to the case that $n_k$ is uniformly bounded by
$N_1$. Let $G:=\max\{-f(n)|n_0\leq n\leq N_1\}$. Suppose $C'>0$. There are
positive integers $T,T'$ with $T'$ depending on $X, P$ and $r$ only, such that
$C'=T/T'$. Choose an integer $N_2>\max\{N_1,G\.T'+N_1\}$.
Assume $N>N_2$. Then for all $n\geq N$
$$h^0(\kf(n))-P(r',\eps,n)\leq-f(n_k-1)-(n+1-n_k)\.C'\leq G-(N_2-N_1)\.C'<0.$$
Again we can restrict to the case $C'=0$. But this gives (1).

Let
$N_3=\lceil N_2+(2g-2)/H^2+1\rceil$ and assume $N>N_3$. Then for all $n\geq
N$, one has
$n>n_k+(2g-2)/H^2$ so that $h^1(\kh_n|_H)=h^1(\kf(n)|_H)=0$. In particular
$$H^2(\kf(n))=H^2(\kf(n+1))=H^2(\kf(n+2))=\ldots,$$
and these cohomology groups must vanish for $n\gg 0$, hence already for $n\geq
N$. This is assertion (2). Moreover,
$$h^0(\kf(n_0)|_H)\leq
h^0(\kf(N_3)|_H)=\chi(\kf(N_3)|_H)=\Delta\chi_{\kf}(N_3)=
\Delta P(r',\eps,N_3)$$
according to (1). Let $M:=\max\{\lceil \Delta P(r',\eps,N_3)\rceil|0\leq
r'\leq r\}$. Then (3) holds.

It remains to prove (4). Since $\kf(N_3)=\kh_{N_3}$, there are short exact
sequences
$$\ses{\ko_H(\nu-N_3)^{\oplus r'}}{\kf(\nu)|_H}{Q_{N_3}}$$
for all $\nu=n_0,\ldots,N_3$. Hence $$h^1(\kf(\nu)|_H)\leq
r'\.h^1(\ko_H(\nu-N_3))\leq r'\.h^1(\ko_H(n_0-N_3))$$
and
$$h^2(\kf(n_0))\leq h^2(\kf(N_3))+\sum_{\nu=n_0+1}^{N_3}h^1(\kf(\nu)|_H)\leq
r(N_3-n_0)\.h^1(\ko_X(n_0-N_3))$$
is uniformly bounded. Since by assumption $h^1(\kf(n_0))\leq Q$ and
$h^0(\kf(n_0))\leq Q$, the Euler characteristic $\chi(\kf(n_0))$
lies in a finite set of integers. By (1) $\Delta\chi_{\kf}$ is given. Hence
$\chi_{\kf}$ lies in a finite set of polynomials. Using (3) and
criterion \ref{boundcrit} we conclude that the set of modules $\kf$ we are
left with is bounded. Therefore there is a constant $N_4>N_3$ such that
$h^1(\kf(n))=0$ if $n\geq N_4$. The lemma holds, if we choose any $N>N_4$.\qed

An immediate consequence of this lemma is the following boundedness result:

\begin{corollary}\label{bound} Suppose $X$ is a surface.
The set of isomorphism classes of $\ko_X$-modules $\ke$ which
occur in $\mu$-semistable pairs $(\ke,\alpha)$ with $\Tors(\kea)=0$
is bounded.
\end{corollary}

\prf Apply lemma \ref{techlem} with $Q=h^0(\ke_0)$. The proof of the lemma
shows that $h^0(\ke(n_0))\leq Q$. By Serre's theorem
$h^0(\ke(n))=\chi_{\ke}(n)=P(r,1,n))$ for all large
enough numbers $n$, so the second alternative of the lemma holds. Part (3)
then states: $h^0(\ke(n_0)|_H)\leq M$
for some $\ke$-regular hyperplane section $H$ and some constant $M$ which is
independent of $\ke$. Therefore the Kleiman criterion applies to the set
of modules $\ke(n_0)$ with the constant $K:=\max\{Q,M,r\.\deg X\}$.\qed

As a consequence of the corollary there is an integer $\hat N$ such that
$\ke(n)$ is globally generated and $h^i(\ke(n))=0$ for all $i>0$, $n\geq\hat N$
and for all $\ko_X$-modules
$\ke$ satisfying the hypotheses of the corollary. Note that according to
lemma $\ref{firstobs}$ among these all
the modules occuring in semistable pairs can be found.

After these preparations we can prove the equivalent to theorem
\ref{boundcurve} in the surface case:

\begin{theorem}\label{boundsurf} Suppose $X$ is a surface. There is an integer
$N$ depending on $X$, $\ko_X(1)$, $h^0(\ke_0)$ and $P$, such that
\begin{itemize}
\item[i)] if $(\ke,\alpha)$ is (semi)stable (with respect to $\delta$) then
$(\ke(n),\alpha(n))$ is sectional (semi)\-stable (with
respect to $\delta(n)$) for all $n\geq N$, and
\item[ii)] if $(\ke(n),\alpha(n))$ is sectional (semi)stable for some $n\geq
N$, then $(\ke,\alpha)$ is (semi)stable.
\end{itemize}
\end{theorem}

\prf By the boundedness result \ref{bound} the dimension of $H^1(\ke(n_0))$ is
uniformly bounded for all $\ke$ satisfying the hypotheses of the corollary.
Let $Q:=h^0(\ke_0)+\max\{h^1(\ke(n_0))\}$. Let
$N$ be the number obtained by applying lemma \ref{techlem}. Without loss of
generality $N>\hat N$.

ad i): Suppose $(\ke,\alpha)$ is (semi)stable. Apply lemma \ref{techlem} to
$\ke$. Since by Serre's theorem $h^0(\ke(n))=\chi(\ke(n))=P(r,1,n)$ for all
sufficiently large $n$, the second alternative of the lemma holds and shows
$h^1(\ke(n))=h^2(\ke(n))=0$ for $n\geq N$. Hence $V:=H^0(\ke(n))$ has
dimension $\chi(n)$. Now let $\kf$ be a submodule of $\ke$. Then either
$h^0(\kf(n))<P(\rk\kf,\eps(\kf),n)$ for all $n\geq N$, in which case we are
done, or we have $\Delta\chi_{\kf}=\Delta(\rk\kf,\eps(\kf))$.
Let $\ke'=\kea$ if $\kf\subseteq\kea$ and $\ke'=\ke$ else.
Let $\ks:=\ke'/\kf$, $\bar\ks:=\ks/\Tors(\ks)$
and let $\bar\kf$ be the kernel of the epimorphism $\ke'\rightarrow\bar\ks$.
Then $\rk\kf=\rk\bar\kf$, $\eps(\kf)=\eps(\bar\kf)$, and we
must have $\Delta\chi_{\bar\kf}=\Delta P(\rk\kf,\eps(\kf))=\Delta\chi_{\kf}$.
Hence $\bar\kf/\kf=\Tors(\ks)$ has zero-dimensional support. There is a short
exact sequence $$\ses{\bar\kf(n_0)}{\ke'(n_0)}{\bar\ks(n_0)}.$$
Now $\bar\ks(n_0)$ cannot have global sections. For otherwise there is a
submodule in $\bar\ks$ isomorphic to $\ko_X(-n_0)$. Let $\kg$ be its preimage
in $\ke'$. Then
$$\Delta\chi_{\kg}=\Delta\chi_{\bar\kf}+\Delta\chi_{\ko_X(-n_0)}\leq
\Delta P(\rk\kf+1,\eps(\kf))=\Delta P(\rk\kf,\eps(\kf))+\Delta P(1,0)$$
contradicting lemma \ref{n0}. But this shows that
\begin{eqnarray*}
h^1(\bar\kf(n_0))\leq h^1(\ke'(n_0))
&\leq&h^1(\ke(n_0))+h^0(\ke(n_0)/\ke'(n_0))\\
&\leq&h^1(\ke(n_0))+h^0(\ke_0(n_0))\leq Q.
\end{eqnarray*}
By part (4) of lemma \ref{techlem} we now conclude that
$$h^0(\kf(n))\leq
h^0(\bar\kf(n))=\chi(\bar\kf(n))\,(\leq)\,P(\rk\kf,\eps(\kf))$$ for all
$n\geq N$ if $\bar\kf\,(\subseteq)\,\ke$. Only the case $\bar\kf=\ke$ for
stable pairs needs special attention: In this case one has
$h^0(\kf(n))<h^0(\ke(n))$, because $\kf$ is a proper submodule of $\ke$ and
$\ke(n)$ is globally generated for all $n\geq N$. Hence
(semi)stability implies sectional (semi)stability for all $n\geq N$.

ad ii) Suppose $(\ke(n),\alpha(n))$ is sectional (semi)stable for some $n\geq
N$. Assume that there exists
a submodule $\kf\subseteq\ke$ with $\Delta\chi_{\kf}>\Delta P(\rk\kf,
\eps(\kf))$.
If such a module exists at all, we may assume that it is maximal with this
property among the submodules of $\ke$. Let
$\ks=\ke/\kf$. The maximality of $\kf$ implies that $\ks$ is torsion free if
$\eps(\kf)=1$ and that $\alpha$ embeds $\Tors(\ks)$ into $\Tors(\ke_0)$ if
$\eps(\kf)=0$. Hence $h^0(\Tors(\ks))\leq Q$. Suppose $\kg$ is any submodule
of $\ks$. Let $\kf'$ be the preimage of $\kg$ under the map $\ke\rightarrow
\ks$. Then $$\Delta\chi_{\kg}+\Delta\chi_{\kf}=\Delta\chi_{\kf'}\leq\Delta
P(\rk\kf+\rk\kg,1)=\Delta P(\rk\kf,\eps(\kf))+\Delta(\rk\kg,1-\eps(\kf)).$$
The inequality in the middle of this line is infered from the maximality of
$\kf$. Hence
$$\Delta\chi_{\kg}\leq\Delta(\rk\kg,1-\eps(\kf))+
\{\Delta P(\rk\kf,\eps(\kf))-\Delta\chi_{\kf}\}
<\Delta P(\rk\kg,1-\eps(\kf)).$$
Therefore we can apply lemma \ref{techlem} to the module $\ks$ with
$\eps=1-\eps(\kf)$. But we did
assume that $(\ke(n),\alpha(n))$ was sectional semistable. Hence there exists
a vector space $V\subseteq H^0(\ke(n))$ of dimension $\chi(n)$ such that
$$\dim(V\cap H^0(\kf(n)))\leq P(\rk\kf,\eps(\kf),n)$$
and
$$h^0(\ks(n))\geq\dim V-\dim(V\cap H^0(\kf(n)))
\geq\chi(n)-P(\rk\kf,\eps(\kf),n)
=P(\rk\ks,1-\eps(\kf),n)$$
This excludes the first alternative of the lemma, and we get
$\Delta\chi_{\ks}=\Delta(\rk\ks,1-\eps(\kf))$ and equivalently
$\Delta\chi_{\kf}=\Delta(\rk\kf,\eps(\kf))$, which contradicts the original
assumption. Thus we have proven that $\Delta\chi_{\kf}\leq\Delta
P(\rk\kf,\eps(\kf))$. But this means that $\ke$ satisfies the
hypotheses of corollary \ref{bound}. By the remark
following the corollary we have $h^0(\ke(\nu))=\chi(\ke(\nu))$ for all
$\nu\geq N$ since $N\geq\hat N$, so that necessarily $V=H^0(\ke(n))$.
Applying lemma \ref{techlem} to $\kf$ we see that either
\begin{itemize}\item[]
$h^0(\kf(\nu))<P(\rk\kf,\eps(\kf),\nu)$ for all $\nu\geq N$, in particular
$\chi_{\kf}<P(\rk\kf,\eps(\kf))$,\end{itemize}
or
\begin{itemize}\item[]
$h^2(\kf(\nu))=0$ for all $\nu\geq N$ and hence
$$\chi_{\kf}(n)=h^0(\kf(n))-h^1(\kf(n))\leq
h^0(\kf(n))\,(\leq)\,P(\rk\kf,\eps(\kf),n),$$
which together with $\Delta\chi_{\kf}=\Delta P(\rk\kf,\eps(\kf))$ implies
$\chi_{\kf}\,(\leq)\,P(\rk\kf,\eps(\kf))$.
\end{itemize}
This finishes the proof.\qed




\subsection{The basic construction}\label{Kapitel13}
Let $X$ be a curve or a surface. By the results of the previous section the
set of modules $\ke$ with fixed Hilbert polynomial $\chi$ that occur
in semistable pairs is bounded. In particular, there is a projective open and
closed part $A$ of the Picard scheme $\Pic(X)$ such that $[\det\ke]\in A$ for
all $\ke$ in semistable pairs. Let $\kl\in\Pic(A\times X)$ be a universal line
bundle. Then there is an integer $N$ such that for all $n\geq N$ the
following conditions are simultaneously satisfied:
\begin{itemize}
\item[-] $0<\delta(n)<\chi(n)$.
\item[-] $\ke$ is globally generated and $h^i(\ke(n))=0$ for all $i>0$
and for all $\ke$ in semistable pairs.
\item[-] $(\ke,\alpha)$ is (semi)stable (with respect to $\delta$) if and only
if $(\ke(n),\alpha(n))$ is sectional (semi)\-stable (with respect to
$\delta(n)$).
\item[-] If $p_A,p_X$ denote the projection maps from $A\times X$ to $A$ and
$X$, respectively, then $R^ip_{A*}(\kl\otimes p_X^*\ko_X(n))=0$ for all
$i>0$,
$\ku_n:=p_{A*}(\kl\otimes p_X^*\ko_X(n))$ is locally free and
$p_A^*(\ku_n)\otimes p_X^*\ko_X(-n)\rightarrow\kl$ is surjective.
\end{itemize}
By twisting the pairs $(\ke,\alpha)$ with $\ko_X(n)$ for sufficiently large
$n$ we can always assume that the assertions above hold for $N=0$. We make
this assumption for the rest of this section and write $p:=\chi(0)$ and
$\bar\delta:=\delta(0)$.

Let $V$ be a vector space of dimension $p$ and let $V_X=V\otimes_k
\ko_X$. Quotient modules of $V_X$ with Hilbert polynomial $\chi$ are
parametrized by a projective scheme $\Quot^{\chi}_{X/V_X}$ (\cite[3.1.]{Gr}).
On
the product $\Quot^{\chi}_{X/V_X}\times X$ there
is a universal quotient $\qtilde:p_X^*V_X\epimorph\ketilde$.
Forming the determinant bundle of $\ketilde$ induces a morphism
$$det:\Quot_{X/V_X}^{\chi}\longrightarrow\Pic(X)$$ so that
$\det\ketilde=(det\times\id_X)^*(\kl)\otimes p_{Quot}^*(\km)$
for some line bundle $\km\in\Pic(\Quot^{\chi}_{X/X_V})$.
Let $Q$ denote the preimage of $A$ under the map $det$. We use the same
symbols for the universal quotient and its restriction to $Q\times X$.

Further let $P:=\bP(\Hom(V,H^0(\ke_0))\dual)$. Again there is a universal
homomorphism $\atilde:(V\otimes_kH^0(\ke_0)\dual)\otimes\ko_P\epimorph
\ko_P(1)$. For sufficiently high $n$ the direct image sheaf
$\kh:=p_{Q*}(\Ker\,\qtilde\otimes p_X^*\ko_X(n))$ is locally free and the
canonical homomorphism
$$\beta:p_Q^*\kh\rightarrow\Ker\,\qtilde\otimes p_X^*\ko_X(n)$$
is surjective, so that there is an exact sequence
$$p_Q^*\kh\otimes
p_X^*\ko_X(-n)\rpfeil{10}{\tilde\beta}p_X^*V_X\rpfeil{10}{\qtilde}\ketilde
\rpfeil{10}{}0.$$
$\tilde\beta$ induces a homomorphism of $\ko_Q$-modules
$$\gamma:\kh\otimes_kH^0(\ke_0(n))\dual\rightarrow\ko_Q\otimes_k(V\otimes_k
H^0(\ke_0)\dual).$$
Let $\ki$ be the ideal in the symmetric algebra
$\ks^*(V\otimes_kH^0(\ke_0)\dual)\otimes_k\ko_Q$ which is generated by the
image of $\gamma$ and let $B\subset P\times Q$ be the corresponding closed
subscheme. Let $\pi_P:B\rightarrow P$ and $\pi_Q:B\rightarrow Q$ be the
projection maps and let $\ko_B(1):=\pi_P^*\ko_P(1)$. This scheme $B$ is the
starting point for the construction of the moduli space for semistable pairs.
We introduce the following notations: Let
$$q_B:=(\pi_Q\times\id_X)^*\qtilde:V\otimes\ko_{B\times
X}\rightarrow\ke_B:=(\pi_Q\times\id_X)^*\ketilde$$
and
$$a_B:=\pi_P^*\atilde:(V\otimes
H^0(\ke_0)\dual)\otimes\ko_B\rightarrow\ko_B(1).$$
By definition of $B$ an arbitrary morphism $h:T\rightarrow P\times Q$ factors
through the closed immersion $B\rightarrow P\times Q$ if and only if the
pull-back under $h$ of the composition $p_P^*\atilde\circ p_Q^*\gamma$ is the
zero map. This is equivalent to saying that the pull-back under $h\times\id_X$
of the induced homomorphism $V\otimes\ko_{P\times Q\times X}\rightarrow
p_P^*\ko_P(1)\otimes p_X^*\ke_0$ factors through $V\times\ko_{T\times
X}\rightarrow(h\times\id_X)^*\ketilde$. This applies in particular to $B$
itself. Let $\alpha_B:\ke_B\rightarrow p_B^*\ko_B(1)\otimes p_X^*\ke_0$ be the
induced homomorphism.

\begin{lemma} (i) There is an open subscheme $Q^0$ of $Q$ such that $u$ is a
point in $Q^0$ if and only if $h^i(\ketilde_u)=0$ for all $i>0$ and the
homomorphism $\tilde q_u:V_X\otimes k(u)\rightarrow\ketilde_u$ induces an
isomorphism on the spaces of global sections.

(ii) Let $(\ke,\alpha)$ be a flat family of pairs parametrized by a
noetherian $k$-scheme $T$. Then there is an open subscheme $S\subseteq T$ such
that $\Ker(\alpha_t)$ is torsionfree for a geometric point $t$ of $T$ if and
only
if $t$ is a point of $S$.
\end{lemma}

\prf (i) By semicontinuity of $h^i$ there is an open subscheme of $Q$ of points
$u$ for which the higher cohomology groups of $\ketilde_u$ vanish. For those
points
$h^0(\ketilde_u)=p$ and hence $H^0(q_u)$ is an isomorphism if and only if
$h^0(\Ker\,q_u)=0$, which again is an open condition for $u$.

(ii) For $n$ large enough there is a locally free $\ko_T$-module $\kg$ and a
surjection\\$\kg\otimes\ko_X(-n)\epimorph\ke\dual$ and dually an inclusion
$\beta':\ke\ddual\rightarrow\kg\dual\otimes\ko_X(n)$. Note that there is an
open subscheme $O$ of $T\times X$ which meets every fibre $X_t$
and for which the restriction $\ke|_O$
is locally free, so that in particular $\vartheta:\ke\rightarrow
\ke\ddual$ is an isomorphism when restricted to $O$. If we let
$\beta=\beta'\circ\vartheta$, then the kernel of $\beta_t:\ke_t\rightarrow
\kg\dual(t)\otimes\ko_X(n)$ is precisely the torsion part of $\ke_t$. Hence
the kernel of $\gamma_t:=(\alpha_t,\beta_t):\ke_t\rightarrow\ke_0\oplus(\kg
\dual(t)\otimes\ko_X(n))$ is the torsion submodule of $\Ker(\alpha_t)$. It
is enough to show that the points $t$ with $\Ker(\gamma_t)=0$ form an open set.
But this is \cite[Cor IV 11.1.2]{EGA}.\qed

Let $S$ be the open subscheme of $B$ which according to the lemma
belongs to the family $(\ke_B,\alpha_B)$, and let $B^0=S\cap(P\times Q^0)$.
The algebraic group $\Sl(V)$ acts naturally on $Q$ and $P$ from the right. On
closed points $[q:V\otimes\ko_X\rightarrow\ke]$ and $[a:V\rightarrow
H^0(\ke_0)]$ this action is given by $[q]\cdot
g=[q\circ(g\otimes\id_{\ko_X})]$ and $[a]\cdot g=[a\circ g]$.

\begin{lemma} $B^0$ is invariant under the diagonal action of $\Sl(V)$ on
$P\times Q$.\end{lemma}

\prf This is clear from the characterization of $B$ as the subscheme of points
$([q],[a])$ for which there is a commuting diagram
$$\begin{array}{ccc}
V\otimes\ko_X&\rpfeil{15}{q}&\ke\\
\spfeil{12}{a}&&\spfeil{12}{\alpha}\\
H^0(\ke_0)\otimes\ko_X&\rpfeil{15}{ev}&\ke_0.
\end{array}$$\qed

$B^0$ has the following local universal property:

\begin{lemma}\label{locunivprop} Suppose $T$ is a noetherian $k$-scheme
parametrizing a
flat family  $(\ke,\alpha)$ of semistable pairs on $X$. Then there is an open
covering $T=\bigcup T_i$ and for each $T_i$ a morphism $h_i:T_i\rightarrow
B^0$
and a nowhere vanishing section $s_i$ in $h_i^*\ko_B(1)$ such that the pair
$(\ke,\alpha)|_{T_i}$ is isomorphic to the pair $((f_i\times\id_X)^*\ke_B,
(f_i\times\id_X)^*(\alpha_B)/s_i)$.
\end{lemma}

\prf Let $T$ be a noetherian scheme and $(\ke,\alpha)$ a flat family of
semistable pairs on $X$ parametrized by $T$. According to the remarks in the
first paragraph of
this section the direct image sheaf $p_{T*}\ke$ is locally free of rank $p$
\cite[Thm 12.8]{h1}. Hence locally on $T$  there are trivializations
$V\otimes\ko_T\rightarrow p_{T*}\ke$, which lead to quotient maps
$q:V\otimes\ko_{T\times X}\rightarrow \ke$. By the universal property of $Q$
there is a $k$-morphism $f:T\rightarrow Q$ and a uniquely determined
isomorphism $\Phi:(f\times\id_X)^*\ketilde\rightarrow\ke$ such that
$\Phi\circ(f\times\id_X)^*\qtilde=q$. Moreover, the composition
$$V\otimes\ko_{T\times X}\rpfeil{10}{q}\ke\rpfeil{10}{\alpha}p_X^*\ke_0$$
determines a homomorphism $a:V\otimes\ko_T\rightarrow H^0(\ke_0)\otimes\ko_T$.
By the universal property of $P$ there is a morphism $g:T\rightarrow P$ and a
uniquely determined nowhere vanishing section $s$ in $g^*\ko_P(1)$ such that
$a=g^*\atilde/s$. It is clear from the construction that
$h:=(f,g):T\rightarrow P\times Q$ factors through $B^0$. $\Phi^{-1}$ is an
isomorphism from $\ke$ to $(f\times\id_X)^*\ketilde=(h\times\id_X)^*\ke_B$,
and $\alpha\circ\Phi=(h\times\id_X)^*(\alpha_B)/s$.\qed

If $h:T\rightarrow B$ and $g:T\rightarrow\Sl(V)$ are morphisms let $h\cdot g$
denote the composition $T\rpfeil{12}{(h,g)}B\times\Sl(V)\rightarrow B$, where
the last map is the induced group action of $\Sl(V)$ on $B$.

\begin{lemma}\label{induce} Suppose $T$ is a noetherian $k$-scheme and
$h=(f,g):T\rightarrow
B^0\subset P\times Q$ a $k$-morphism. $h$ induces (locally) isomorphism
classes of families of pairs. If $g:T\rightarrow\Sl(V)$ is
a morphism, then the families induced by $h$ and $h\cdot g$ are isomorphic.
Conversely, if $h_1$ and $h_2$ induce isomorphic families parametrized by $T$,
then there is an etale morphism $c:T'\rightarrow T$ and a
morphism $g':T'\rightarrow\Sl(V)$ such that the
morphisms $(h_1\circ c)\cdot g'$ and $(h_2\circ c)$ are equal.
\end{lemma}

\prf Let $h:T\rightarrow B^0$ be a $k$-morphism. Applying $(h\times\id_X)^*$
to $\ke_B$ and $\alpha_B$ induces a family $\ke_T$ and a homomorphism
$\alpha_T:\ke_T\rightarrow h^*(\ko_B(1))\otimes p_X^*\ke_0$. Locally there are
nowhere vanishing sections in $h^*\ko_B(1)$. Dividing $\alpha_T$ by any
of these sections defines families of pairs. Two such sections differ by a
section in $\ko_T^*$. But this sheaf embeds into the sheaf of automorphisms of
$\ke_T$. Hence the families induced by different sections are isomorphic. The
second statement is clear. For the third assume that $h_1$ and $h_2$ are
morphisms
such that for $i=1,2$ there are nowhere vanishing sections $s_i\in H^0(T,h_i^*
\ko_B(1))$. Let
$$\ke_i:=(h_i\times\id_X)^*\ke_B,
\quad q_i:=(h_i\times\id_X)^*q_B\quad\hbox{and}\quad
\alpha_i:=(h_i\times\id_X)^*\alpha_B/s_i.$$
Assume that there is an
isomorphism $\Phi:(\ke_1,\alpha_1)\rightarrow(\ke_2,\alpha_2)$ of pairs.
The quotient maps $q_i$ induce isomorphisms $\bar q_i:V\otimes\ko_T\rightarrow
p_{T*}\ke_i$ because of the definition of $B^0$ (\cite[Thm 12.11]{h1}).
The composition $\bar q_2^{-1}\circ p_{T*}\Phi\circ \bar q_1$ corresponds to a
morphism $g:T\rightarrow \Gl(V)$. Define morphisms $c$ and $\ell$ by the
fibre product diagram
$$\begin{array}{ccc}
T'&\rpfeil{20}{c}&T\\
\spfeil{}{\ell}&&\spfeil{}{\det(g)}\\
\Gm&\rpfeil{20}{p^{th}\,\,power}&\Gm
\end{array}$$
and let $g':=(g\circ c)/\ell:T'\rightarrow\Sl(V)$.
It is easy to check that $(h_1\circ c)\.g'=(h_2\circ c)$.
\qed

$\qtilde$ induces a homomorphism $\Lambda^r(\ko_Q\otimes
V_X)\rightarrow\det\ketilde=(det\times\id)^*(\kl)\otimes p_Q^*\km$ and hence
a homomorphism $\Lambda^rV\otimes_k\km\dual\rightarrow
det^*\ku_0=det^*p_{A*}\kl$ (\cite{Ma}). This finally leads to morphisms
$T:Q\rightarrow P':=\bP(\kh om(\Lambda^rV,\ku_0)\dual)$ and
$\tau:=(\pi_P,T):B\rightarrow P\times P'$.

\begin{lemma} $\Sl(V)$ acts naturally on $P'$ from the right, $T$ and $\tau$
are equivariant morphisms with respect to this action.\qed
\end{lemma}

We can choose a very ample line bundle $\kn$ on $A$ such that
$\kn':=\ko_{P'}(1)\otimes p_A^*\kn$ is
very ample on $P'$. For any positive numbers $\nu,\nu'$ the
line bundle $\ko_P(\nu)\otimes(\kn')^{\otimes\nu'}$ is very ample on
$P\times P'$ and inherits a canonical linearization with respect to the
$\Sl(V)$-action \cite[1.4,1.6]{m10}. Choose $\nu$ and $\nu'$such that
$\nu/\nu'=r\bar\delta/(p-\bar\delta)$.
Let $Z^{(s)s}\subseteq P\times P'$ be the open subscheme of (semi)stable
points with respect to this linearization. Here \em stable \em means \em
properly stable \em in the sense of Mumford.

\begin{theorem}\label{mainresult} The open subscheme
$B^{(s)s}=B^0\cap\tau^{-1}(Z^{(s)s})$ of $B$
has the following property: A morphism $h:T\rightarrow B^0$ induces families
of
(semi)stable pairs in the sense of lemma \ref{induce} if and only if $h$
factors through $B^{(s)s}$. The restriction of $\tau$
to $B^{ss}$ is a finite morphism $\taus:B^{ss}\rightarrow Z^{ss}$.
\end{theorem}

For the proof we need a stability criterion for $\tau([a],[q])$, and we
need it in slightly greater generality. But before this, note that if
$q:V_X\rightarrow\ke$ defines a point $[q]$ in $Q(k)$, then the fibre of the
projective bundle $P'$ through the point $T([q])$ is isomorphic to
$P'':=\bP(Hom(\Lambda^rV,H^0(\det\ke))\dual)$, and $\tau([a],[q])$ is a
(semi)stable point in $P\times P'$ if and only if it is (semi)stable point in
$P\times P''$ with respect to the canonical linearization of
$\ko_P(\nu)\otimes\ko_{P''}(\nu')$ (\cite[4.12]{Ma}). In particular, the
choice of $\kn$ is of no
consequence for the definition of $Z^{(s)s}$.

\begin{proposition}\label{prop} Let $(\ke,\alpha)$ be a pair with $\det\ke
\in A$ and torsionfree $\kea$. Suppose there is a generically surjective
homomorphism $q:V_X\rightarrow\ke$ such that $q\circ\alpha\neq 0$. Let
$T:\Lambda^rV\rightarrow H^0(\det\ke)$
and $a:V\rightarrow H^0(\ke_0)$ be the derived homomorphisms. Then $\Ta$ is a
(semi)stable point in $P\times P''$ with respect to the given linearization
if and only if $q$ injects $V$ into $H^0(\ke)$ and $(\ke,\alpha)$ is sectional
(semi)stable with respect to $\bar\delta$.
\end{proposition}

The proof of this proposition is postponed to the next section.

\em Proof of theorem \ref{mainresult}: \em Pairs $(\ke,\alpha)$ that
correspond to points $([a],[q])$ in $B^0$ satisfy the hypotheses of
the proposition, for $q$ is surjective, $H^0q$ isomorphic and $\kea$
torsionfree. Hence by proposition \ref{prop} and theorems \ref{boundcurve} and
\ref{boundsurf} $(\ke,\alpha)$ is (semi)stable if and only if $\tau([a],[q])$
is a (semi)stable point. This proves the
first assertion of the theorem. In order to show that $\taus:=\tau|_{B^{ss}}$
is a finite morphism it is enough to show that $\taus$ is proper and injective
 (\cite[IV 8.11.1]{EGA}). This will be done in two steps:

\begin{proposition}\label{properness} $\taus$ is a proper morphism.
\end{proposition}

\prf Using the valuation criterion it suffices to show the following:
Let $\bar C=Spec\,R$ be a nonsingular affine curve, $c_0\in\bar C$ a closed
point defined by a local parameter $t\in R$ and $C$ the open complement of
$c_0$. Suppose we are given a commutative diagram
$$\begin{array}{ccc}C&\rpfeil{15}{h}&B^{ss}\\
\spfeil{}{\iota}&&\spfeil{}{\tau^{ss}}\\
\bar C&\rpfeil{15}{m}&Z^{ss}.
\end{array}$$
We must show that (at least locally near $c_0$) there is a lift $\bar h:\bar
C\rightarrow B^{ss}$ of $m$ extending $\bar h$. Making $C$ smaller if
necessary we may assume that $h$ induces homomorphisms
$$\ko_C\otimes V_X\rpfeil{12}{q}\kf\rpfeil{12}{\alpha}\ko_C\otimes\ke_0,$$
so that $(\kf,\alpha)$ is a flat family of semistable pairs. Using Serre's
theorem one can find a locally free $\ko_X$-module $\kh$ and an epimorphism
$\ko_C\otimes\kh\dual\epimorph\kf\dual$. The kernel of the dual homomorphism
$\beta:\kf\rightarrow\ko_C\otimes\kh$ is the torsion submodule $\Tors(\kf)$.
Since $\Ker\,\alpha$ and $\Ima\,\alpha$ are $C$-flat, $(\Ker\,\alpha)_c\subset
(\Ker\,\alpha_c)$.
Since the kernel of the restriction of $\alpha$
to any fibre $X\times c$, $c\in C$, is torsion free by lemma \ref{firstobs},
$\Ker\,\alpha$ is also torsion free. Therefore
$$(\alpha,\beta):\kf\rightarrow\ko_C\otimes(\ke_0\oplus\kh)$$
is injective. There are integers $a,b$ such that the composition
$$\ko_C\otimes
V_X\rpfeil{15}{q}\kf\rpfeil{20}{(t^a\alpha,t^b\beta)}\ko_C\otimes
(\ke_0\oplus\kh)$$
extends to a homomorphism
$$\lambda:\ko_{\bar C}\otimes V_X\longrightarrow\ko_{\bar C}\otimes(\ke_0\oplus
\kh)$$ which is nontrivial in each component when restricted to the special
fibre $X_{c_0}$. Let $\bar\kf$ be the maximal submodule of
$\ko_{\bar C}\otimes(\ke_0\oplus\kh)$ with the properties
\begin{center}$\bar\kf|_{C\times X}=\kf$, ${\rm Im}\,\lambda\subseteq
\bar\kf$ and $\dim\Supp(\bar\kf/{\rm Im}\,\lambda)<e$;
\end{center}
and let $\bar\alpha:\bar\kf\rightarrow\ko_{\bar C}\otimes\ke_0$ be the
projection map. Then $\bar\kf$ is $\bar C$-flat,
$(\bar\kf,\bar\alpha)|_{C\times X}\isom(\kf,\alpha)$ and $q_{c_0}:V_X
\rightarrow\kf_{c_0}$ is generically surjective. Moreover $\bar\alpha_{c_0}$
is nonzero and $\Ker\,\bar\alpha_{c_0}$ is torsion free. For assume that
$\Tors(\Ker\,\bar\alpha_{c_0})\neq0$ and let $\tilde\kf$ be the kernel of the
composite epimorphism
$$\bar\kf\epimorph\kf_{c_0}\epimorph\kf_{c_0}/\Tors(\Ker\,\bar\alpha_{c_0}).$$
Then there is a short exact sequence
$$\ses{\bar\kf}{\tilde\kf}{\Tors(\Ker\,\bar\alpha_{c_0})}.$$
By construction $\bar\alpha$ extends to $\tilde\alpha:\tilde\kf\rightarrow
\ko_{\bar C}\otimes\ke_0$. Since $\kh$ is normal and the codimension
of $\Supp\,\Tors(\Ker\,\alpha)$ in $\bar C\times X$ is greater than 1,
$\bar\beta$ also
extends to a homomorphism $\tilde\beta:\tilde\kf\rightarrow\ko_{\bar C}\otimes
\kh$. Finally $(\tilde\alpha,\tilde\beta):\tilde\kf\rightarrow\ko_{\bar
C}\otimes(\ke_0\oplus\kh)$ is injective, contradicting the maximality of
$\bar\kf$. Hence indeed $\Tors(\Ker\,\bar\alpha_{c_0})=0$. Since $q_{c_0}$
is generically surjective, $\Ker\,\bar\alpha_{c_0}$ torsionfree and $\bar
\alpha_{c_0}\circ q_{c_0}\neq 0$, we can apply proposition \ref{prop} to the
pair
$(\bar\kf_{c_0},\bar\alpha_{c_0})$. By
assumption on the map $m$ the induced point in $P\times P'$ is semistable,
hence
$H^0q_{c_0}$ is injective and $(\kf_0,\alpha_0)$ sectional
semistable. But then necessarily $\kf_{c_0}$ is globally generated,
$H^0q_{c_0}$
isomorphic and $q_{c_0}$ surjective. This shows that $h$ extends to a
morphism $\bar h:\bar C\rightarrow B$ with $\bar h(c_0)\in B^{ss}$.
\qed

\begin{proposition}\label{monomorphism} $\taus$ is injective.
\end{proposition}

\prf Assume that for $i=1,2$ there are closed points $([a_i:V\to H^0(\ke_0)],
[q_i:V_X\to\ke_i])$  with the same image under $\tau$.
We may assume that $a_1=a_2$ and $\det\,\ke_1=\det\,\ke_2$. Then there is an
open subscheme
 $\emptyset\not=U\subset X$ such that ${\ke_1}_{|U}$, ${\ke_2}_{|U}$ are
locally free and
are in fact isomorphic as quotients of ${V_X}_{|U}$.
Then $\ke_1/\Tors(\ke_1)$ and
$\ke_2/\Tors(\ke_2)$ are isomorphic as quotients of $V_X$ via a map
$\Phi:\ke_1/\Tors(
\ke_1)\to\ke_2/\Tors(\ke_2)$ (\cite{Ma}, lemma 4.8).
The kernels of the induced homomorphisms $\alpha_i:\ke_i\to\ke_0$ are
torsionfree,
so that the natural map $\ke_i\to\ke_0\oplus\ke_i/\Tors(\ke_i)$ are
injective. The diagram
$$\begin{array}{ccccc}
V_X&\epimorph&\ke_1&\lra&\ke_0\oplus\ke_1/\Tors(\ke_1)\\
\|&&&&\spfeil{1}{\id+\Phi}\\
V_X&\epimorph&\ke_2&\lra&\ke_0\oplus\ke_2/\Tors(\ke_2)
\end{array}$$
commutes and shows $\ke_1$ and $\ke_2$ are isomorphic as quotients of
$V_X$.\qed

This completes the proof of theorem \ref{mainresult} up to the proof of
proposition \ref{prop}.\qed

\begin{theorem}\label{mainthm}
Assume that $X$ is a smooth projective  variety of dimension one or two.
Then there is a
projective $k$-scheme $\km_\delta^{ss}(\chi,\ke_0)$ and a natural
transformation
$$\varphi:\underline\km_\delta^{ss}(\chi,\ke_0)\lra\Hom_{\Spec\, k}(
{}~~~,\km_\delta
^{ss}(\chi,\ke_0))~~,$$
such that $\varphi$ is surjective on rational points and
$\km_\delta^{ss}(\chi,\ke_0)$
is minimal with this property. Moreover, there is an open
subscheme $\km_\delta^s
(\chi,\ke_0)\subset\km_\delta^{ss}(\chi,\ke_0)$ such that $\varphi$ induces
an isomorphism of
subfunctors$$\underline\km_\delta^s(\chi,\ke_0)\stackrel{\cong}{\lra}
\Hom_{\Spec\, k}(~~~~,\km_\delta^s(\chi,\ke_0))~~,$$
i.e. $\km_\delta^s(\chi,\ke_0)$ is a fine moduli space
for all stable pairs.\end{theorem}

\prf By \cite[1.10]{m10} and \cite{g2} there is a projective $k$-scheme
$\km^{ss}$
and a morphism $\rho:B^{ss}\lra\km^{ss}$ which is a good quotient for the
$\Sl(V)$-action
on $B^{ss}$. By lemma \ref{induce} and theorem \ref{mainresult} any family
of semistable pairs parametrized by $T$ induces morphisms $T_i\to B^{ss}$
for an
appropriate open covering $T=\bigcup T_i$ such that the composition with
$\rho$ glue
to a well-defined morphism $T\to\km^{ss}$. This establishes a natural
transformation
$$\varphi:\underline\km_\delta^{ss}(\chi,\ke_0)\lra\Hom_{\Spec\, k}(~
{}~,\km^{ss})~~.$$
If $\psi:\underline\km^{ss}\lra\Hom(~~,N)$ is a similar transformation, then
the family
$(\ke_B,\alpha_B)|_{B^{ss}}$ induces an $\Sl(V)$-invariant morphism
$B^{ss}\lra N$, which
therefore factors through a morphism $\km^{ss}\lra N$. Moreover there is
an open subscheme
$\km^s\subset\km^{ss}$ such that $B^s=\rho^{-1}(\km^s)$ and
$\rho|_{B^{s}}:B^s\lra\km^s$
is a geometric quotient. In order to see that the family
$(\ke_B,\alpha_B)|_{B^s}$ descends
to give a universal pair on $\km^s$ it is enough to show that the stable
pairs have no
automorphism besides the identity.
But assume that $\Phi\not=\id$ is an automorphism of a stable pair
$(\ke,\alpha)$,
i.e. $\Phi:\ke\stackrel{\cong}{\lra}\ke$ and $\alpha\circ\Phi=\alpha$.
Then $\psi=\Phi-\id$ is a nontrivial homomorphism from $\ke$ to $\kea$.
Apply the stability conditions to $\Ker\psi\subset\ke$ and
$\Ima\psi\subset\kea$
to get$$\rk\,\ke\cdot\chi_{\Ker\psi}<\rk(\Ker\psi)(\chi_\ke-\delta)+\delta
\cdot\rk\,\ke$$
and
$$\rk\,\ke\cdot\chi_{\Ima\psi}<\rk(\Ima\psi)(\chi_\ke-\delta)~.$$
Summing up and using $\chi_{\Ima\psi}+\chi_{\Ker\psi}=\chi_\ke$ and
$\rk(\Ima\psi)
+\rk(\Ker\psi)=\rk\ke$ we get the contradiction $\chi_\ke<\chi_\ke$.\qed


\subsection{Geometric stability conditions}\label{Kapitel14}
In this section we prove proposition \ref{prop}.
Let $q:V_X\rightarrow\ke$ and $\alpha:\ke\rightarrow\ke_0$ be homomorphisms of
$\ko_X$-modules. To these data we can associate vector space homomorphisms
$T:\Lambda^rV\rightarrow H^0(\det\,\ke)$ and $a:V\rightarrow H^0(\ke_0)$.
If $q$ is generically surjective, then $T$ is nontrivial, and if
$\alpha\circ q\neq 0$, then $a$ is nontrivial. Let $\Ta$
denote the corresponding closed point in $P\times P''$ (notations as in
section \ref{Kapitel13}).

The group $\Sl(V)$ acts on $P\times P''$ by
$$\Ta\cdot g=([a\circ g], [T\circ\Lambda^rg]).$$
We want to investigate the stability properties of $\Ta$ with respect to
an $\Sl(V)$-linearization of the very ample line bundle
$\ko_{P\times P''}(\nu,\nu')$, where $\nu,\nu'$ are positive integers.
These stability properties depend on the ratio
$\eta:=\nu/\nu'$ only. We will make use of the Hilbert criterion to
decide about (semi)\-stability.
Let $\lambda:\Gm\rightarrow\Sl(V)$ be a 1-parameter subgroup, i.\ e.\ a
nontrivial group homomorphism. There is a basis $v_1,\ldots, v_p$ of $V$ such
that $\Gm$ acts on $V$ via $\lambda$ with weights
$\gamma_1,\ldots,\gamma_p\in\zz$:
$$\lambda(u)\.v_i=u^{\gamma_i}\.v_i\qquad\hbox{for all }u\in\Gm(k).$$
Reordering the $v_i$ if necessary we may assume that
$\gamma_1\leq\ldots\leq\gamma_p$, $\sum\gamma_i=0$, since ${\rm
det}\lambda=1$, and $\gamma_1<\gamma_p$, since $\lambda\neq 1$.

For any multiindex $I=(i_1,\ldots,i_r)$ with $1\leq i_1<\ldots<i_r\leq p$ let
$v_I=v_{i_1}\wedge\ldots\wedge v_{i_r}$ and
$\gamma_I=\gamma_{i_1}+\cdots+\gamma_{i_r}$. The vectors $v_I$ form a basis of
$\Lambda^rV$, and $\Sl(V)$ acts with weights $\gamma_I$ with respect to this
basis. $T(v_I)\neq0$ if and only if the sections
$q(v_{i_1}),\ldots,q(v_{i_r})$ are generically linearly independent, i.\ e.\
generate $\ke$ generically. Now let
$$\mu=\mu([a],\lambda):=-{\rm min}\{\gamma_i|a(v_i)\neq 0\}.$$
$$\mu'=\mu([T],\lambda):=-{\rm min}\{\gamma_I|T(v_I)\neq 0\}$$

\begin{lemma}[Hilbert criterion] $\Ta$ is a (semi)stable point in
$P\times P''$ with respect to $\ko(\nu,\nu')$ if and only if
$\hat\mu:=\eta\.\mu+\mu'(\geq)0$ for all 1-parameter subgroups $\lambda$.
\end{lemma}

\prf \cite[Thm 2.1.]{m10}\qed

For any linear subspace $W\subset V$ let $\ke_{(W)}\subset\ke$ be the
submodule which is characterized by the properties :
$\ke/\ke_{(W)}$ is torsionfree and $\ke_{(W)}$ is generically generated by
$q(W\otimes\ko_X)$. In particular, let $\ke_{(i)}=\ke_{(\langle
v_1,\ldots,v_i\rangle)}$, $i=0,\ldots,p$ for a given basis $v_1,\ldots,v_p$.
Then there is a filtration
$$\Tors(\ke)=\ke_{(0)}\subset\ke_{(1)}\subset\ldots\subset\ke_{(p-1)}\subset
\ke_{(p)}=\ke.$$
Since $\ke_{(i)}/\ke_{(i-1)}$ is torsionfree, one has either
$\ke_{(i)}=\ke_{(i-1)}$ or $\rk\ke_{(i)}>\rk\ke_{(i-1)}$. Consequently, there
are integers $1\leq k_1<\ldots<k_r\leq p$ marking the points, where the rank
jumps, i.\ e.\ $k_{\rho}$ is minimal with $\rk\ke_{(k_{\rho})}=\rho$. Let $K$
denote the multiindex $(k_1,\ldots,k_r)$. If $I$ is any multiindex as above,
let $i_0=0$ and $i_{r+1}=p+1$ for notational convenience.

\begin{lemma} $\mu'=-\gamma_K$.\end{lemma}

\prf By construction $T(v_K)\neq 0$. We must show that $\gamma_K\leq\gamma_I$
for every multiindex $I$ with $T(v_I)\neq 0$. For any $I$ and any
$t\in\{1,\ldots,r\}$ we let $\ke_{I,t}=\ke_{\langle
v_{i_1},\ldots,v_{i_t}\rangle)}$. Now suppose $T(v_I)\neq 0$. Let
$\ell={\rm max}\{\lambda|k_t=i_t\quad\forall t<\lambda\}$. If $\ell\geq
r+1$, then $I=K$ and we are done. We will procede by descending induction on
$\ell$. By definition of $K$, we have $k_{\ell}<i_{\ell}$. Define
$\ke_{I,t}'=\ke_{(\langle v_{k_1},\ldots,v_{k_{\ell}},v_{i_{\ell}},\ldots
v_{i_t}\rangle)}$ for $t=\ell,\ldots p$. Then $\ke_{I,t}\subset\ke'_{I,t}$,
and $t\leq\rk\,\ke'_{I,t}\leq t+1$. Let $m={\rm min}\{t|\rk\,\ke'_{I,t}=t,\ell
\leq t\leq p\}$. Now define a multiindex
$$I'=(k_1,\ldots,k_{\ell},i_{\ell},\ldots,i_{m-1},i_{m+1},\ldots,i_p).$$
(If $m=\ell$, drop the $i_{\ell},\ldots,i_{m-1}$ part; if $m=p$, drop the
$i_{m+1},\ldots,i_p$ part.) Then we have $T(v_{I'})\neq 0$, and
$\gamma_{I'}\leq\gamma_I$ by monotony of $I$ and $\gamma$. Moreover, $I'$ and
$K$ agree at least in the first $\ell$ entries. Thus by induction
$\gamma_K\leq\gamma_{I'}\leq\gamma_I$.\qed

Let $\ell:=\min\{i|a(v_i)\neq0\}$. Obviously $\mu=-\gamma_{\ell}$, so that
$\hat\mu=-\gamma_K-\eta\.\gamma_{\ell}$. Now $\ell$ and $K$ depend on the basis
$v_1,\ldots,v_p$ only, and $\mu$ is a linear function of $\gamma$ for fixed
$\ell$ and $K$. Using these notations, the Hilbert criterion can be expressed
as follows:

\begin{lemma} $\Ta$ is a (semi)stable point if and only if
$$\min_{\hbox{\em\scriptsize bases of
$V$}}\min_{\gamma}-(\gamma_K+\eta\.\gamma_{\ell})\quad
(\geq)\quad0.$$\qed
\end{lemma}

We begin with minimizing over the set of all weight vectors $\gamma$. This is
the cone spanned by the special weight vectors
$$\gamma^{(i)}=(\underbrace{i-p,\ldots,i-p}_{i},\underbrace{i,\ldots,i}_
{p-i})$$
for $i=1,\ldots,p-1$. For any weight vector $\gamma$ can be expressed as
$\gamma=\sum_{i=1}^{p-1}c_i\gamma^{(i)}$ with nonnegative rational coefficients
$c_i=(\gamma_{i+1}-\gamma_i)/p$. In order to check (semi)stability for a given
point it is enough to show $\hat\mu(\geq)0$ for each of these basis vectors.
Let $\delta_i=1$ or $0$ if $\ell\leq i$ or $>i$,
respectively. Evaluating $\hat\mu$ on $\gamma^{(i)}$ we get numbers
$$\mu^{(i)}=p\.(\max\{j|k_j\leq i\}+\eta\.\delta_i)-i\.(r+\eta).$$
If $i$ increases, $\mu^{(i)}$ decreases unless $i$ equals $\ell$ or any of the
numbers $k_j$, in which case $\mu^{(i)}$ might jump. The critical values of
$i$ therefore are $\ell-1$ and $k_j-1$, $j=1,\ldots,r$, and the corresponding
critical values of $\mu^{(i)}$ are:

\begin{tabular}{ll}
$p\.(j-1)-(k_j-1)\.(r+\eta)$& if $1\leq j\leq r$, $1<k_j\leq\ell$,\\
$p\.(j-1)-(\ell-1)\.(r+\eta)$& if $1\leq j\leq r+1$, $k_{j-1}<\ell\leq
k_j$, $1<\ell$,\\
$p\.(j-1+\eta)-(k_j-1)\.(r+\eta)$& if $1\leq j\leq r$, $\ell<k_j$.
\end{tabular}

If we put $\ell_j=\min\{k_j,\ell\}$, then the conditions imposed by these
values of $\hat\mu$ can be comprised as follows:

\begin{tabular}{lll}
(1)& $0\,(\leq)\,p\.(j-1)-(\ell_j-1)\.(r+\eta)$ &if $1\leq
j\leq r+1$, $1<\ell_j$\\
(2)& $0\,(\leq)\,p\.(j-1+\eta)-(k_j-1)\.(r+\eta)$ &if $1\leq j\leq r$.
\end{tabular}

In the next step one should minimize these terms over all bases of $V$. But in
fact, the relevant information is not the used basis itself but the flag of
subspaces of $V$ which it generates. The stability criterion takes the
following form:

\begin{lemma}\label{1form} $\Ta$ is a (semi)stable point if and only if
\begin{itemize}
\item[1)] $\dim W\.(r+\eta)\,(\leq)\,p\.\rk\,\ke_{(W)}$ for all subspaces
$0\neq
W\subseteq\Ker\,a$.
\item[2)] $\dim W\.(r+\eta)\,(\leq)\,p\.(\rk\,\ke_{(W)}+\eta)$ for all
subspaces
$0\neq W\subseteq V$ with $\rk\,\ke_{(W)}\,(\leq)\,r$.
\end{itemize}\qed\end{lemma}

We give the stability criterion still another form, shifting our attention
from subspaces of $V$ to submodules of $\ke$:

\begin{lemma}\label{2form} $\Ta$ is a (semi)stable point if and only if
\begin{itemize}
\item[(0)] $H^0q$ is an injective map.
\item[(1)] $V\cap H^0\kf=0$ or $\dim(V\cap H^0\kf)\.(r+\eta)\,(\leq)\,
p\.\rk\kf$ for all submodules $\kf\subseteq\Ker\alpha$.
\item[(2)] $\dim(V\cap H^0\kf)\.(r+\eta)\,(\leq)\,p\.(\rk\kf+\eta)$ for all
submodules $\kf\subseteq\ke$ with \newline $\rk\kf\,(\leq)\,\rk\ke$.
\end{itemize}\end{lemma}

\prf If $\Ta$ is semistable, let $W:=\Ker H^0q$. Then $W\subseteq\Ker\,a$.
{}From the lemma above it follows that $\dim W\leq
p/(r+\eta)\.\rk\,\ke_{(W)}=0$.
Hence (0) is a necessary condition. It is to show that the conditions(1) and
(2) of lemma \ref{1form} and of lemma \ref{2form} are equivalent.
Suppose we are given a submodule
$\kf\subseteq\ke$. Let $W:=V\cap H^0\kf$. Then $q(W\otimes\ko_X)\subseteq\kf$
and $\rk\,\kf=\rk\,\ke_{(W)}$. Moreover, if $\kf\subseteq\kea$, then
$W\subseteq\Ker\,a$. Now either $W=0$ or \ref{1form} applies and gives
\ref{2form}. Conversely,
if $W\subseteq V$ is given, let $\kf:=q(W\otimes\ko_X)$. Then $W\subseteq V\cap
H^0\kf$ and $\rk\,\ke_{(W)}=\rk\,\kf$. Again, if $W\subseteq\Ker\,a$, then
$\kf\subseteq\kea$. Hence \ref{2form} implies \ref{1form}.

Finally, we replace $\eta$ by a more suitable parameter:
$$\bar\delta=\frac{p\.\eta}{r+\eta}\qquad\eta=\frac{r\.\bar\delta}
{p-\bar\delta}.$$
Since $\eta$ was a positive rational number, $\bar\delta$ is
confined to the open interval $(0,p)$, which of course tallies with the
data of the previous section. The following theorem, which differs from
proposition \ref{prop} only in the choice of words, summarizes the discussion
of this section:

\begin{theorem} If in addition to the global assumptions of this section
$\kea$
is torsionfree, then $\Ta$ is a (semi)stable point of $P\times P''$
if and only if the following conditions are satisfied:
\begin{itemize}
\item[-] $H^0q$ is an injective homomorphism.
\item[-] $(\ke,\alpha)$ is sectional stable with respect to $\bar\delta$.
\end{itemize}\end{theorem}

\prf If $\kea$ is torsionfree then every nontrivial submodule of $\kea$ has
positive rank. Hence condition (1) in \ref{2form} can be replaced by
\begin{itemize}\item[(1')] $dim(V\cap H^0\kf)\.(r+\eta)(\leq)p\.\rk\,\kf$ for
all
submodules $\kf\subseteq\kea$.
\end{itemize}
As a result of replacing $\eta$ by $\bar\delta$ in (1') and \ref{2form}(0),(2)
one obtains the defintion of sectional (semi)stability.\qed


\section{Applications}
This chapter is organized as follows. In \ref{Kapitel21} we show that the
existence of semistable pairs
gives an upper bound for
$\delta$.  Rationality conditions on $\delta$ imply the
equivalence of semistability and stability. If $\delta$ varies within certain
regions the semistability
conditions remain unchanged. This is formulated and
specified  for the rank two case.

\ref{Kapitel22} deals with Higgs pairs. Again we concentrate on the rank two
case. We make
the first
step to generalize the diagrams of Bertram and Thaddeus to algebraic
surfaces. The
restriction of $\mu-$stable \vbs on an algebraic surface to a curve of
high degree
induces an immersion of the \ms  of \vbs on the surface into the \ms  of
\vbs on the curve.
The understanding of this process is important, e.g. for the computation of
Donaldson polynomials
and for the study of the geometry of the moduli space on the surface
(\cite{T}). With the help of a restriction
theorem for $\mu-$stable pairs $(\ke,\alpha:\ke\to\ko)$ we construct an
approximation
of this immersion, which will hopefully shed some light
on the relation between the original moduli spaces. It is remarkable that
the limit of any approximation is
independent of the polarization.

In \ref{Kapitel23} we first compare our \st  for $\ke_0=\ko_D^{\oplus r}$,
where $D$ is a divisor on a curve,
with the notion of Seshadri of stable sheaves with level structure along a
divisor(\cite{Se}).
We will have a closer look at the \ms  of rank two sheaves of degree 0 with
a level structure
at a single point. Furthermore certain results from \ref{Kapitel22} are
reconsidered in the case of $\ke_0$
being a vector
bundle on a divisor.


\subsection{Numerical properties of $\delta$}\label{Kapitel21}
Let $X$ be a smooth projective variety with an ample divisor $H$, $\ke_0$ a
coherent $\ko_X-$module and
$\delta$ a positive rational polynomial of degree $\dim X-1$ with leading
coefficient $\delta_1\geq0$.

\begin{lemma}\label{rest}Assume $(\ke,\alpha)$ is a semistable pair such that
$\kea\not=0$. Then

$$\delta(\leq)\chi_\ke-\frac{\rk\ke}{\rk\kea}(\chi_\ke-\chi_{\ke_0})~~.$$

If $\ke_0\cong\ko_X$ and $\rk\ke>1$, then $$\delta(\leq)\frac{\rk\ke\.\chi_
{\ko_X}-\chi_\ke}{\rk\ke-1}$$
and in particular $$\delta_{1}(\leq)-\frac{\deg\ke}{\rk\ke-1}~~.$$

If $\ke_0$ is torsion, then $$\delta(\leq)\chi_{\ke_0}$$
and in particular  $\delta_1(\leq)\deg{\ke_0}$.\end{lemma}
\prf The first inequality follows immediately from the stability
condition i). If $\ke_0\cong\ko_X$ use
$\chi_{\kea}=\chi_\ke-\chi_{\Ima\alpha}\geq\chi_\ke-\chi_{\ke_0}$ and
$\rk\kea=\rk\ke-1$.\qed

It is much more convenient to work with $\mu-$\st only. In fact for the
general $\delta$ one can
achieve that every semistable pair is $\mu-$stable.

\begin{lemma}\label{ratio} There exists a discrete set of rationals
$0\leq...<\eta_i<\eta_{i+1}
<...$ including $0$, such that  for $\delta_1\in (\eta_i,\eta_{i+1}) $
every semistable pair
with
respect to $\delta$ is in fact $\mu-$stable and the $\mu-$stability
conditions depend only on
$i$.\end{lemma}
\prf Define $\{\eta_i\}:=[0,-{d}/({r-1}))\cap\{({ar-sd})
/({r-s})|a,s\in\IZ,~0\leq s<r\}$.
If $\delta_1\in(\eta_i,\eta_{i+1})$,
then the right hand sides of the $\mu-$ semistability conditions
$\deg\kg\leq{sd}/{r}-
\delta_1{s}/{r}$ and $\deg\kg\leq{sd}/{r}+\delta_1({r-s})/{r}$ are
not integer ($s= \rk\kg$). Therefore $\mu-$semistability and $\mu-$stability
coincide. Moreover, the integral
parts
of the right hand sides depend only on $i$, i.e. for two different choices
of $\delta_1$ in the
intervall $(\eta_i,\eta_{i+1})$ the $\mu-$stability conditions are the
same.\qed

More explicit results can be achieved in special cases:
\begin{proposition}\label{indforlb}For $r=2$ and $\ke_0\in Pic(X)$ and
$\delta_1\in(\eta_i,\eta_i+2)$,where $\eta_i:=
\max\{0,2i+d\}$ with $i\in\IZ$, every semistable pair is $\mu-$stable.
The stability in this
region does not depend on $\delta$.\end{proposition}
\prf For $\ke_0\in Pic(X)$ all semistable pairs $(\ke,\alpha)$ have
torsionfree $\ke$ and
$\rk\kea=1$. In particular the \st conditions concern rank one subsheaves
only. Now $\delta_1
\in(\eta_i,\eta_i+2)$ is equivalent to $-1-i<{d}/{2}-{\delta_1}/{2}<-i$,
$i+d-1<{d}/
{2}+{\delta_1}/{2}<i+d$ and $\delta_1>0$.\qed

As the last numerical criterion we mention
\begin{lemma} Assume $\delta_1<\min_{0\leq s<r}\{(r-sd)/({r-s})
+{r}({r-s})[{sd}/{r}]\}$.
\begin{itemize}
\item[i)] Then every sheaf $\ke$ in a semistable pair $(\ke,\alpha)$ without
torsion in dimension zero
is torsionfree and $\mu-$semistable.
\item[ii)] If $\ke$ is torsionfree and $\mu-$semistable and
$\alpha:\ke\to\ke_0$ a nontrivial
homomorphism such that $\kea$ does not contain a destabilizing
subsheaf, then $(\ke,\alpha)$ is $\mu-$stable.
\end{itemize}\end{lemma}
\prf The condition on $\delta_1$ is equivalent to either of the two conditions:

$[sd/r,sd/r+\delta_1(r-s)/r)\cap\IZ=\emptyset$ for $0\leq s<r$.

$[sd/r-\delta_1/r,sd/r)\cap\IZ=\emptyset$ for $0<s\leq r$.\qed

\subsection{Higgs pairs in dimension one and two}\label{Kapitel22}
A Higgs pair in this context is a \vb together with a global section.
(This notion should not be confused
 with a Higgs field as a section $\theta\in H^0(\ke nd\ke\otimes\Omega^1_X)$
with $\theta\wedge\theta=0$!)
Instead of considering a global section we prefer to work with a homomorphism
 from the dualized bundle to the structure sheaf. These objects  will be
called pairs as in the
general context.

First we remind of the situation in the curve case, which was motivation
for us to go on.

\begin{definition} Let $C$ be a smooth curve. As introduced in \ref{Kapitel13}
 $\km_\delta^{ss}(d,2,\ko)$(resp. $\km_\delta^{ss}(\kq,2,\ko)$) denotes the
moduli space of
semistable pairs $(\ke,\alpha:\ke
\to\ko)$ with respect to $\delta$,
where $\ke$ is a rank two sheaf of degree $d$ (with
determinant $\kq$).\end{definition}

\begin{remark}\label{torsfree}Notice, that $\delta$ is just a number and that
a sheaf occuring in a semistable
pair is always torsionfree and hence a vector bundle. Moreover the stability
conditions reduce to
$\deg(\kea)\leq{d}/{2}-{\delta}/{2}$ and $\deg(\kg)\leq{d}/{2}+{\delta}/{2}$
for all line bundles $\kg\subset\ke$.\end{remark}

For the following we assume $d<0$.

\begin{definition}\label{UC}$U_{C,i}(d):=\km_\delta^{ss}(d,2,\ko)$
and
$SU_{C,i}(\kq):=\km_\delta^{ss}(\kq,2,\ko)$, where $\delta\in(
max\{0,2i+d\},2i+d+2)$.\end{definition}

Note that according to proposition \ref{indforlb} the spaces
$U_{C,i}(d)$ and $SU_{C,i}(\kq)$ do not
depend on the choice of $\delta$
\begin{proposition}(M. Thaddeus) $U_{C,i}(d)$ and $SU_{C,i}(\kq)$ are
projective fine
moduli spaces. Every semistable pair is automatically
stable.\end{proposition}
\prf \cite{Th} or \ref{mainthm}\qed

\begin{proposition}\begin{itemize}
\item[i)] For $i\geq-d$ the \mss $U_{C,i}(d)$ are empty.
\item[ii)] For $i=\lfloor-{d}/{2}-1\rfloor+1$ there are
morphisms $$U_{C,i}(d)\lra U(d)$$
and$$SU_{C,i}(\kq)\lra SU(\kq)~~,$$where $U(d)$ and $SU(\kq)$ are
the \mss of semistable \vbs
of degree $d$ and determinante $\kq$, resp. The fibre over a stable
bundle $\ke$
is isomorphic to $\IP(H^0(\ke\dual)\dual)$. In particular they are
projective bundles for $0\gg d\equiv1(2)$.
\item[iii)] A pair $(\ke,\alpha)$ lies in $SU_{C,-d-1}(\kq)$ if and
only if there is a nonsplitting
exact sequence of the form
$$\sesq{
\kq}{}{\ke}{\alpha}{\ko}~.$$Thus $SU_{C,-d-1}\cong\IP({{\rm Ext}}^1(\ko,
\kq)\dual)$.\end{itemize}
\end{proposition}
\prf i) and ii) follow from the general criteria. A similar result as iii)
holds in the surface
case. We give the proof there.\qed

The following picture illustrates the situation:$$\begin{array}{clll}
SU_{C,\lfloor-{d}/{2}-1\rfloor+1}(\kq)&SU_{C,\lfloor-{d}/{2}-1\rfloor+2}(
\kq)&....
&SU_{C,-d-1}(\kq)\cong\IP({{\rm Ext}}^1(\ko,\kq)\dual)\\
\downarrow&&&\\
SU(\kq)&&&
\end{array}$$

M. Thaddeus is able 'to resolve the picture' by a sequence of blowing ups
and downs.
In particular all the spaces $SU_{C,i}$ are rational. This process makes
it possible to trace
a generalized theta divisor on $SU_{C,i}$ to a certain section of $\ko(k)$
on $\IP(H
^1(\kq))$. This method is used in \cite{Th} to give a proof of the Verlinde
formula.

We go on to proceed in a similar way in the case of a surface.

Let $X$ be an algebraic surface with an ample divisor $H$.
Now $\km_\delta^{ss}(d,c_2,2,\ko)$ ($
\km_\delta^{ss}(\kq,c_2,2,\ko)$) denotes the moduli space of semistable pairs
$(\ke,\alpha:\ke\to\ko_X)$ with respect to $\delta$, where $\ke$ is a rank
two sheaf of degree $d$ ($:=c_1.H$) (with
determinant $\kq$) and second Chern class $c_2$. For the existence of such
pairs it is necessary
that $\delta$ be linear with nonnegative leading coefficient $\delta_1$.
 As in \ref{torsfree} a sheaf occuring in a semistable pair is
torsionfree and the stability conditions are
$$\chi_\kg(\leq)\frac{\chi_\ke}{2}-\frac{\delta}{2}$$
for all rank one subsheaves $\kg\subset\kea$ and
$$\chi_\kg(\leq)\frac{\chi_\ke}{2}+\frac{\delta}{2}
$$for all rank one subsheaves $\kg\subset\ke$.

\begin{definition}\label{UCX}For $\delta$ such that $\delta_1
\in(\max\{0,2i+d\},2i+d+2)$ we define $U_i:=\km_\delta^{ss}(d,c_2,2,\ko)$ and
$SU_i:=\km_\delta^{ss}(\kq,c_2,2,\ko)$.\end{definition}
Again, note that according to \ref{indforlb} the definition does not
depend on the choice of $\delta$.

\begin{corollary} $U_i$ and $SU_i$ are projective fine moduli spaces.
Every semistable pair is $\mu-$stable.\end{corollary}
\prf It follows immediately fom \ref{mainthm} and section
\ref{Kapitel21}.\qed

\begin{proposition}If $(\ke,\alpha)$ is a $\mu-$semistable pair with
respect to
$\delta$, then $4c_2(\ke)-c_1^2(\ke)\geq-{\delta_1}/({4H^2})$.\end{proposition}
\prf If $(\ke,\alpha)$ is a $\mu-$semistable pair the homomorphism
$\alpha$ can be extended
to a homomorphism $\ke^{{\rm vv}}\to \ko$ and the resulting pair is
still $\mu-$semistable with
$c_1(\ke^{{\rm vv}})=c_1(\ke)$ and $c_2(\ke^{{\rm vv}})\leq c_2(\ke)$. Thus
it is enough to prove the
inequality for  locally free pairs. If $\ke$ itself is a $\mu-$semistable
bundle the Bogomolov
inequality says $4c_2-c_1^2\geq0$. If $\ke$ is not $\mu-$semistable,
then there is
an exact sequence $$\ses{\kl_1}{\ke}{\kl_2\otimes I_{Z}}~~,$$where $I_Z$
is the ideal sheaf of a
zero dimensional subscheme and $\kl_1$ and $\kl_2$ are line bundles with
${\deg\ke}/{2}<
\deg\kl_1\leq{\deg\ke}/{2}+({1}/{2})\delta_1$ and ${\deg\ke}/{2}-({1}/{2})
\delta_1\leq\deg\kl_2<{\deg\ke}/{2}$. Using
$c_2(\ke)=c_1(\kl_1)c_1(\kl_2)+l(Z)\geq c_1(
\kl_1)c_1(\kl_2)=({1}/{4})\{(c_1(\kl_1)+c_1(\kl_2))^2-(c_1(\kl_1
)-c_1(\kl_2))^2\}=({1}/{4})
c_1^2(\ke)-\frac{1}{4}(c_1(\kl_1)-c_1(\kl_2))^2$ and Hodge index theorem,
which gives
$(c_1(\kl_1)-c_1(\kl_2))^2\leq ({(\deg\kl_1-\deg\kl_2)^2})/{H^2}$ we infer
the claimed
inequality. Notice, that for $\delta_1\to0$ the inequality converges to
the usual
Bogomolov inequality.\qed

\begin{proposition}\begin{itemize}
\item[i)] For $i\geq-d$ the \mss $U_i$ and $SU_i$ are empty.
\item[ii)]
If $i=\lfloor-{d}/{2}-1\rfloor+1$, then every pair $(\ke,\alpha)\in U_i$ has
a $\mu-$semistable
$\ke$. There is rational map $U_i\to U(c_1,c_2)$ (the \ms of semistable,
torsionfree sheaves),
which is a morphism for $d\equiv1(2)$. The image of the rational map
contains all $\mu-$stable sheaves
$\ke$ with ${{\rm Hom(\ke,\ko)}}\not=0$. The fibre over such a point is
$\IP({{\rm Hom}}(\ke,\ko)\dual)$.
\item[iii)] Every pair $(\ke,\alpha)\in SU_{-d-1}$ sits in an nontrivial
extension of the form
$$\sesq{I_{Z_1}\otimes\kq}{}{\ke}{\alpha}{I_{Z_2}}~~,$$where $I_{Z_i}$
are the ideal sheaves of
certain zero dimensional subscheme. In the case $Z_1=\emptyset$, e.g.
$\ke$ is locally free, every
such extension gives in turn a stable pair $(\ke,\alpha)\in
SU_{-d-1}$.\end{itemize}\end{proposition}

\prf {\it i)} and {\it ii)} follow again from \ref{Kapitel21}
If $(\ke,\alpha)\in SU_{-d-1}$, then $\deg\kea< d+{1}/{2}$, which is
equivalent to $
\deg(\Ima\alpha)>-{1}/{2}$. Since $\Ima\alpha\subset\ko$ it follows
$\Ima \alpha= I_{Z_2}$.
 A splitting of the induced exact sequence would lead to the
contradiction $0\leq\deg I_{Z_2}
\leq-{1}/{2}$. Let $(\ke,\alpha)$ be given by a sequence with $
Z_1=\emptyset$. For $\kg\subset\kea$ one gets the required inequality
$\deg\kg\leq d<d+{1}/{2}$.
If $\kg\subset\ke$ and $\kg\not\subset\kea$, the sheaf $\kg$ has the
form $\kg= I_{Z_3}\subset I_{Z_2}$.
Without restriction we can assume that $\ke/\kg$ is torsionfree.
Since $\ke/\kg$ is an extension of $
I_{Z_2}/I_{Z_3}$ by $\kq$ and ${{\rm Ext}}^1(I_{Z_2}/I_{Z_3},\kq)=0$,
$\kg$ in fact equals $I_{Z_2}$
and therefore defines a splitting of the sequence.\qed

\begin{corollary}The set of all pairs $(\ke,\alpha)\in SU_{-d-1}$
with $\kea$ locally free, which
in particular contains all  locally free pairs, forms a projective scheme
over $Hilb^{c_2}(X)$
with fibre over $[Z]\in Hilb^{c_2}(X)$ isomorphic to $\IP({{\rm Ext}}^1
(I_Z,\kq)\dual)$. \end{corollary}
\prf If $(\ke,\alpha)$ is a universal family over $SU_{-d-1}\times X$,
then the set of points
$t\in SU_{-d-1}$ with $l((\coker\,\alpha)_t)$ maximal is closed. It is easy
to see that $
(\coker\,\alpha)_t\cong\coker(\alpha_t)$ and that $l(\coker(\alpha_t))$ is
maximal, i.e. is
equal to $c_2$
 if $
\Ker(\alpha_t)$ is locally free. Therefore the set of all pairs with locally
free kernel $\kea$
 is closed and $\ko/\Ima\alpha$ induces the claimed morphism to
$Hilb^{c_2}(X)$.\qed
\begin{corollary}\label{indepofpol} The moduli space of all locally free
pairs $(\ke,\alpha)\in SU_{-d-1}$
does not depend on the polarization of $X$.\end{corollary}

\begin{remark}{\it i)} Bradlow introduced in \cite{Br} the notion of
$\phi-$\st with
respect to a parameter $\tau$. If we set $\delta_1=-d+({\tau}/{2\pi})vol(X)$
($d$ is the degree of $\ke$) both notions
coincide, i.e. a pair $(\ke,\alpha:\ke\to\ko)$ with a locally free $\ke$
is $\mu-$stable in
our sense if and only if $(\ke\dual,\phi=\alpha\dual\in H^0(\ke\dual))$
is $\phi-$stable with respect to the
parameter $\tau$ in Bradlow's sense. He proves a Kobayashi-Hitchin
correspondence in this
situation, i.e. he shows: $(\ke,\alpha)$ is $\mu-$stable ( or a sum of a
$\mu-$stable pair with $\mu-$stable
bundles) if and only if the vortex
equation has a solution, i.e. there exists a hermitian metric $H$ on
$\ke\dual$, such that
$$\Lambda_\omega F_H+\tau\frac{i}{2}id=\frac{i}{2}\phi\otimes\phi^{*_H}~~.$$
$F_H$ is the curvature of the metric connection on $\ke\dual$, $\omega$ is
a fixed K\"ahler form and $\Lambda_\omega$ is the adjoint of $\wedge\omega$.
Now,if $(\ke,
\alpha)\in SU_{-d-1}$ one can take $\delta$ near to $-d$. That corresponds
to $\tau\to0$.
Although \ref{indepofpol} shows that $SU_{-d-1}$ is independent of the
polarization $H$,
i.e. of the Hodge metric, for us there is no obvious reason in the
analytical equation.\\
{\it ii)} In \cite{Rei} the space $SU_{-d-1}$ is stratified and equipped
with certain line
bundles. These objects Reider calls Jacobians of rank two alluding to a
Torelli kind
theorem for algebraic surfaces.
\end{remark}

In order to study the restriction of $\mu-$stable \vbs to curves of
high degree it could
be usefull to study the restriction of $\mu-$stable pairs to those
curves. As a generalization
of a result of Bogomolov we prove

\begin{theorem}\label{restofpairs}
For fixed $c_1,c_2,\delta$ and $H$ there exists a constant $n_0$, such
that for $n\geq n_0$ and
any smooth curve $C\in|nH|$ the
restriction of every locally free, $\mu-$stable pair to $C$ is a
$\mu-$stable pair on the curve
with respect to $n\delta_1$.\end{theorem}
\prf If $\ke$ is locally free the kernel $\kea$ is a line bundle. In
particular the restriction
 of the injection $\kea\subset\ke$ to a curve remains injective.
Thus $(\kea)_C=\Ker(\alpha_C)$.
Since $\deg(\kea)_C=n\deg\kea$, the two inequalities
$\deg\kea<{\deg\ke}/{2}-{\delta_1}/{2}$
and $\deg({\kea})_C<{\deg\ke_C}/{2}-{n\delta_1}/{2}$ are equivalent.
Thus the first of the stability conditions on $C$ is always satisfied. In
order to prove
the second we proceed in two steps.\\
i) By Bogomolov's result (\cite{Bo}) there is a constant $n_0$, such that
the restriction of
a $\mu-$stable \vb to a smooth curve $C\in|nH|$ for $n\geq n_0$ is stable.
Since the inequality
$\deg\kg<{\deg\ke_C}/{2}+{n\delta_1}/{2}$ for a line bundle
$\kg\subset\ke_C$ is
weaker than the stability condition on $\ke_C$, the theorem follows
immediately from Bogomolov's
result for all $\mu-$stable pairs $(\ke,\alpha)$, where $\ke$ is a
$\mu-$stable vector bundle.\\
ii) Therefore it remains to prove the theorem for pairs with $\ke$
not $\mu-$stable. Any
such \vb is an extension of $\kl_2\otimes I_Z$ by $\kl_1$,
where $\kl_1$ and $\kl_2$ are
line bundles with ${\deg\ke}/{2}\leq\deg\kl_1<{\deg\ke}/{2}+({1}/{2})
\delta_1$. $I_Z$
is as usual the ideal sheaf of a zero dimensional subscheme. If $C\in|nH|$
is a curve with
$C\cap Z=\emptyset$, then the restriction of the extension to $C$ induces
the exact sequence
$$\ses{(\kl_1)_C}{\ke}{(\kl_2)_C}~~.$$ If $\kg\subset\ke_C$ is a line bundle,
then either
$\kg\subset(\kl_1)_C$ or $\kg\subset(\kl_2)_C$. This implies
$\deg\kg\leq\deg(\kl_1)_C
=n\deg\kl_1<{\deg\ke_C}/{2}+({1}/{2})n\delta_1$ or $\deg\kg\leq\deg(\kl_2)_C=
n\deg\ke-n\deg\kl_1\leq{\deg\ke_C}/{2}$. Hence $(\ke_C,\alpha_C)$ is stable.
If $C\cap Z\not=\emptyset$ we only get a sequence of the
form $$\ses{(\kl_1)_C(Z.C)}{\ke_C}
{(\kl_2)_C(-Z.C)}$$Notice, that
$\ko_C(-Z.C)\cong(I_Z\otimes\ko_C)/\Tors(I_Z\otimes\ko_C)$. As above
$\deg\kg\leq\deg(\kl_1)_C+\deg(Z.C)\leq\deg(\kl_1)_C+l(Z)$
or $\deg\kg\leq{ \deg(\ke_C)}/{2}
$ for every line bundle $\kg\subset\ke_C$. If
$\deg\kl_1+{l(Z)}/{n}<{\deg\ke}/{2}+
{\delta_1}/{2}$, then $(\ke_C,\alpha_C)$ is stable.
There exists a positive number $\varepsilon$ depending only on the degree,
$\delta$ and $H$,
such that $\deg\kl_1\leq{\deg\ke}/{2}+{\delta_1}/{2}-\varepsilon$. Thus it
suffices to
bound $l(Z)$ by $n_0\varepsilon$. That is done by the following computation.
$l(Z)=c_2-c_1(\kl_1)
c_1(\kl_2)=c_2-{c_1^2}/{4}+({1}/{4})(c_1(\kl_1)-c_1(\kl_2))^2\leq c_2-{c_1^2}/
{4}
+{\delta_1^2}/({4H^2})$. Thus $n_0>({1}/{\varepsilon})(c_2-{c_1^2}/{4}+
{\delta_1^2}/({4H^2}))$ satisfies $l(Z)<n_0\varepsilon$.\qed

With the notation of \ref{UC} and \ref{UCX} one proves
\begin{corollary} For fixed $c_1,c_2$ and $H$ there exists a number $n_0$,
such that for every
smooth curve $C\in|nH|$ for $n_0\leq n\equiv1(2)$ and every
$i$ with $\lfloor{-d}/{2}-1\rfloor+1
\leq i\leq-d-1$ the restriction of pairs gives an injective immersion,
i.e. an injective
morphism with injective tangent map:$$U_i^f\to
U_{C,in+({n-1})/{2}}$$\end{corollary}
(The superscript denotes
the subset of all locally free pairs)\\
\prf The technical problem here is, that the constant $n_0$
in the last theorem depends on
$\delta$ and not only on $i$. Therefore we fix for every $i$
a very special $\delta$, namely
$\delta_1=2i+d+1$. Since we only consider finitely many $i$'s
there is
an $n_0$, such that the restriction gives a morphism
$U_i^f\to U_{C,in+({n-1})/{2}}$. Here we use $n\equiv1(2)$. Since the
occuring family of \vbs is bounded one can choose $n_0$,
such that $H^k(X,\kh om(\ke,\ke')(-nH))=0$ ($k=0,1$) and
$H^0(\ke\dual(-nH))=0$ for $n\geq n_0$ and
all \vbs $\ke$ and $\ke'$ occuring in a pair in one of the \mss $U_i$.
Thus $(\ke,\alpha)_C
\cong(\ke',\alpha')_C$ if and only if $\ke\cong\ke'$ and $\alpha$ maps to
$\alpha'$ under this isomorphism,
i.e. the restriction morphism is injective. A standard argument in
deformation theory shows
that the Zariski tangent space of $U_i^f$ at $(\ke,\alpha)$ is isomorphic
to the hypercohomology
$\IH^1(\ke nd\ke\dual\to\ke\dual)$ of the indicated complex which is given
by $\varphi\mapsto\varphi
(\alpha\dual)$ (\cite{We}). Analogously, the Zariski tangent space of
$U_{C,j}$ at $(\ke_C,\alpha_C)$ is isomorphic
to the hypercohomology $\IH^1(\ke nd\ke_C\dual\to\ke_C\dual)$. The Zariski
tangent map is described
by the restriction of hypercohomology classes. Both hypercohomology groups
sit in
exact sequences of the form$$...\to
H^0(\ke\dual)\to\IH^1(\ke nd\ke\dual\to\ke\dual)\to H^1(\ke nd\ke)\to...$$
and$$...\to
H^0(\ke_C\dual)\to\IH^1(\ke nd\ke_C\dual\to\ke_C\dual)\to
H^1(\ke nd\ke_C)\to...~~~,$$resp. By our
assumptions the restrictions $H^0(\ke\dual)\to H^0(\ke_C\dual)$
and $H^1(\ke nd\ke\dual)\to H^1(\ke nd\ke_C)
$ are injective. Hence the Zariski tangent map of the restriction
of stable pairs is injective,
too.\qed

We remark that neither the starting nor the end point of the series
of \mss on the surface is
sent to the corresponding point of the series \mss on the curve. A slight
generalization of the theorem allows
to restrict $\mu-$stable pairs to a stable pair on a curve $C\in|nH|$  with
respect to
the parameter $n\delta_1+c$, where
$c$ is a constant depending only on $\delta_1,c_1,c_2$ and $H$.


\subsection{Framed bundles and level structures}\label{Kapitel23}
In this paragraph we consider pairs of rank $r$,
where $\ke_0\cong\ko_D^{\oplus r}$ or
more generally where $\ke_0$ is a \vb of rank $r$ on a divisor $D$.\\
We start with pairs on a curve. In this case $D$ is a finite
sum of points. As
far as we know,
Seshadri was the first to consider and to construct \mss for such pairs.
In \cite{Se}
they were called sheaves with a level structure. The general stability
conditions as
developped in this paper and specialized
to this case present a slight generalization of Seshadris stability concept
in terms
of the parameter $\delta$, which in \cite{Se} is always $l(D)$. The geometric
invariant theory
which Seshadri used
to construct the \mss differs from the one in \ref{Kapitel13}. In \cite{Se}
a
point
$[\ko^{\oplus N}\epimorph\ke]$ of the Quotscheme is mapped to a point
$$([\ko^{\oplus N}(x_1)\epimorph
\ke(x_1)],...,
[\ko^{\oplus N}(x_n)\epimorph\ke(x_n)])$$ in the product of Grassmannians
(the $x_i$ are
sufficiently many generic points. The conditions for a point in this product
to be semistable
in the sense of geometric invariant theory translate into the semistability
properties for  pair.
 However, to generalize the construction to
the higher dimensional case one has to map the Quotscheme into a different
projective space as in \ref{Kapitel13} and
study the stability conditions there.

\begin{lemma}If the genus of the curve is at least 2, there exists a
semistable pair of rank $r$
and degree $d$ with respect to $(\ko_D^{\oplus r},\delta)$ if and only if
$0<\delta\leq r\cdot l(D)$.\end{lemma}

\prf The 'only if' part was proven in \ref{rest}, since $r\cdot
l(D)=h^0(\ke_0)$. For the
'if' direction we pick a stable \vb $\ke$ of rank $r$ and degree $d$ and an
isomorphism
$\alpha:\ke_D\cong\ko_D^{\oplus r}$. The induced pair is semistable.\qed

\begin{corollary} The \mss $\km_\delta^{ss}(d,r,\ko_D^{\oplus r})$ of
semistable pairs with
$0<\delta\leq r\cdot l(D)$ exist as projective schemes of generic dimension
$r^2(g-1)+r^2\cdot l(D)$
\end{corollary}(cp. \cite{Se}, III.5., there is a misprint in the dimension
formula in \cite{Se})\\
There are two new features in the theory of pairs compared with the \mss of
vector bundles.
First, to compactify one really has to use sheaves with torsion supported on
$D$. Secondly,
the set of semistable pairs which are not stable may have only codimension 2,
whereas the set of
semistable \vbs which are not stable
is at least $2g-3$ codimensional in the \ms of all semistable vector bundles.
To give an example we  describe the moduli space $\km_1^{ss}(0,2,k(P)^{
\oplus 2})$
 of sheaves of rank two and degree zero with
a level structure at a reduced point $P\in X$ with $\delta=1$. Here we try to
compute the
S-equivalence
in geometric terms, which is not clear to us in the general context.\\
The stability conditions say\begin{itemize}
\item[i)] $\deg\kg(\leq)-\frac{1}{2}$ for all rank one subsheaves
$\kg\subset\kea=\Ker\alpha$.
\item[ii)] $\deg\kg(\leq)\frac{1}{2}$ for all rank one subsheaves
$\kf\subset\ke$.
\item[iii)] $l(\ke/\kea)(\geq)1$
\item[iv)] $l(\Tors(\ke))(\leq)1$
\item[v)] $\alpha$ is injective on the torsion $\Tors(\ke)$.\end{itemize}
Therefore the sheaves $\ke$ occuring in  semistable pairs in
$\km_1^{ss}(0,2,k(P)^{\oplus 2})$ are
either locally free or of the form $\kf\oplus k(P)$ with $\kf$ locally free.\\
First we classify all pairs $(\ke,\alpha)$ with locally free $\ke$.
By ii) such
a bundle $\ke$ has to be semistable as a bundle. If $\ke$ is a stable
bundle, then every
pair $(\ke,\alpha)$ with an arbitrary $\alpha\not=0$ is semistable and is
stable if and only if $\rk(\alpha)=2$,
i.e. $\alpha(P)$ is bijective. If $\ke$ is only semistable there are two
cases to consider:
Either a) $\ke\cong\kl_1\oplus\kl_2$, where $\kl_1$ and $\kl_2$ are line
bundles of degree $0$ or
b) $\ke$ is given as a nontrivial extension of two such line bundles.\\
a) If $\kl_1\cong\kl_2$ then $(\ke,\alpha)$ is semistable if and only if
$\alpha$ is bijective. If
$\kl_1\not\cong\kl_2$ then $(\ke,\alpha)$ is semistable if and only if
none of the
restrictions $\alpha_{|\kl_i
(P)}$  is trivial. \\
If $\kl_1\cong\kl_2$ and $\alpha$ bijective, the pair $(\ke,\alpha)$ is in
fact stable,
since $l(\ke/\kea)=2>1$. If $\alpha$ is only of rank one we can always find
an inclusion
$\kl_1\subset\kl_1\oplus\kl_1$ with $\kl_1=\Ker(\alpha_{|\kl_1})$, which
contradicts i).
For $\kl_1\not\cong\kl_2$ one has to consider line bundles
$\kl\subset\kl_1\oplus\kl_2$ of degree
zero with $\kl=\Ker(\alpha_{|\kl})$, because that is the only possibility to
contradict  i).
But such a line bundle has to be isomorphic to one of the summands with the
natural inclusion.
Therefore the stability condition is equivalent to $\alpha_{|\kl_i}\not=0$.\\
b) If $\ke$ is a nonsplitting extension $$\ses{\kl_1}{\ke}{\kl_2}$$ a pair
$(\ke,\alpha)$ is
semistable iff $\kl_1\not=\Ker(\alpha_{|\kl_1})$. That is, since every line
bundle $\kl\subset\ke$
of degree $0$ either defines a splitting of the sequence or maps
isomorphically to $\kl_1$.\\
The next
step is to determine all semistable pairs $(\kf\oplus k(P),\alpha)$.
Here we claim, that
such a pair is semistable iff $\kf$ is stable and $\alpha_{|k(P)}$ is
injective. Let
$(\kf\oplus k(P),\alpha)$ be semistable and $\kl\subset\kf$ a line
bundle. Then, since
$\kl\oplus k(P)\subset\kf\oplus k(P)$, the semistability conditions
for the pair give
$\deg(\kl\oplus k(P))\leq{1}/{2}$, i.e. $\deg\kl\leq-{1}/{2}={\deg\kf}/{2}$.
Let now
$\kf$ be a stable bundle. If $\kl$ is a rank one subsheaf of
$\kf\oplus k(P)$, then it either
 injectively injects into $\kf$
or has torsion part $k(P)$ and therefore satisfies the
required inequality.\\
Next we look at the isomorphism classes of stable pairs. If $\ke$ is a
stable bundle,
two pairs $(\ke,\alpha)$ and $(\ke,\alpha')$ are isomorphic if and only
if $\alpha$ and $\alpha'$ differ by
a scalar. For $\ke$ of the form $\kl_1\oplus\kl_2$ the automorphism group
of $\ke$ is either
$\IC^*\times\IC^*$ for $\kl_1\not\cong\kl_2$ or $\GL(2)$ for
$\kl_1\cong\kl_2$. In the first case
the set of isomorphism classes of stable pairs for fixed $\ke$ is
isomorphic to
$\PGL(2)/\{{\beta\,0\choose0\,
\gamma}|\beta,\gamma\in\IC^*\}$.
 In the latter case
all stable pairs are isomorphic for fixed $\ke$, they all define the same
point in the
moduli space. If $\ke$
is given by a nonsplitting
exact sequence$$\ses{\kl_1}{\ke}{\kl_2}$$ the automorphism group is either
$\IC^*$ for $\kl_1\not
\cong
\kl_2$ or $\{{\beta\,\gamma\choose0\,\beta}|\beta\in\IC^*,\gamma\in\IC\}$
for $\kl_1\cong\kl_2$.
Therefore every such extension induces either a $\PGL(2)-$family of stable
pairs in the \ms or
a
$\PGL(2)/\{{\beta\,
\gamma\choose0\,\beta}|\beta\in\IC^*,\gamma\in\IC\}-$family of stable pairs in
the moduli space.\\
In order to describe the S-equivalence
we claim, that the orbit of a pair $(\ke,\alpha)$ is closed if and only if
either the pair is stable, i.e. $\ke$ is a semistable \vb and $\alpha$ of
rank two, or $\ke$ is of
the form $\kf\oplus k(P)$ with a stable \vb $\kf$ of degree $-1$ and
$\alpha_{|\kf}=0$.\\
If $\ke$ is locally free and $\alpha$ of rank one there is an extension of
the form
$$\ses{\kea}{\ke}{k(P)}~~.$$If $\psi\in{\rm Ext}^1(k(P),\kea)$ denotes the
extension class
one can easily construct a family of pairs over $\IC\cdot\psi$, which gives
the pair $(\ke,\alpha)$ outside $0$ and $(\kea\oplus k(P),\alpha
\cdot{\rm pr}_{k(P)})$ on the special fibre,
where ${\rm pr}_{k(P)}$ is the projection to $k(P)$. Obviously this pair is
again
semistable. If $(\kf\oplus k(P),\alpha)$ is a semistable pair with
$\alpha=(\alpha_1,\alpha_2)$,
the pair $(\kf\oplus k(P),(t\cdot\alpha_1,\alpha_2))$ converges
constantly to a pair with $
\alpha_{|\kf}=0$ for $t\to 0$.
In order to prove the claim it is therefore enough to show that the orbit of
such a pair is
closed. If there were a family parametrized by a curve with a point $O$,
which outside $O$ were isomorphic to a fixed semistable pair
$(\kf\oplus k(P),\alpha)$ with $\alpha_{|\kf}=0$
and over this point $O$ isomorphic to another pair of this kind,
the family of the kernels would give a family of stable bundles, which would
be constant for all points
except $O$. Since the stable bundles are separated, it has to be constant
everywhere. Finally, using the constance of the images of the maps $\alpha$
outside the point $O$
one concludes that the family of pairs is constant.\\
If $\km_1^s(0,2,k(P)^{\oplus 2})$ denotes the subset of all stable pairs we
summarize the results
in the following proposition

\begin{proposition}\begin{itemize}
\item[i)] $\km_1^{ss}(0,2,k(P)^{\oplus 2})\setminus
\km_1^s(0,2,k(P)^{\oplus 2})\cong \IP_1\times U(-1,2)$,
where $U(-1,2)$ is the \ms of stable rank two \vbs of degree $-1$.
\item[ii)] There is a morphism $\km_1^s(0,2,k(P)^{\oplus 2})\to U(0,2)$,
which is a $\PGL(2)
-$fibre bundle over $U(0,2)^s$ and whose fibre over a point
$[\kl_1\oplus\kl_2]\in U(0,2)\setminus
U(0,2)^s$ is isomorphic to $$\PGL(2)/\{{\beta\,0\choose0\,\gamma}|
\beta,\gamma\in\IC^*\}
\cup\{\PGL(2)\times\IP({\rm Ext}^1(\kl_1,\kl_2)\dual)\}$$
for $\kl_1\not\cong\kl_2$ and isomorphic to
$$\{pt\}\cup\{\PGL(2)/\{{\beta\,\gamma\choose0\,\beta}
|\beta\in\IC^*,\gamma\in\IC^*\}\times\IP({\rm Ext}^1(\kl_1,\kl_2)\dual)\}$$
for $\kl_1\cong\kl_2$.\end{itemize}\end{proposition}
\prf The isomorphism in {\it i)} is given by $(\kf\oplus k(P),\alpha)\mapsto
(\alpha(k(P)),\kf)$.
The morphism in {\it ii)} is induced by the universality property of the
moduli space.\qed

In particular the dimension of $\km_1^{ss}(0,2,k(P)^{\oplus2})$ is $4g$
($g$ is the genus of the curve)
and the dimension of
$\km_1^{ss}(0,2,k(P)^{\oplus2})\setminus\km_1^s(0,2,k(P)^{\oplus2})$ is
$4g-2$. Thus the codimension is two,
independently of the genus.

Finally we want to study the situation in the two dimensional case. Let
$X$ be a surface with
an effective divisor $C$ and $\ke_0$ be a \vb of rank $r$ on $C$. A
framing of a \vb $\ke$ of
rank $r$ on $X$ along $C$ in the strong sense as introduced in
\cite{l4} is an isomorphism
$\alpha:\ke_C\cong\ke_0$. In \cite{l4} the question of the existence of
\mss for such pairs
$(\ke,\alpha)$
was asked ($\alpha$ denotes the isomorphism as well as the composition of
this isomorphism
with the surjection $\ke\epimorph\ke_C$). In fact, under additional
conditions, fine \mss
for such framed bundles were constructed as algebraic spaces. These
additional conditions are:
$C$ is good and $\ke_0$ is simplifying. If $C=\sum b_i C_i$ with prime
divisors $C_i$
and $b_i>0$ $C$ is called good if there exist nonnegative integers $a_i$,
such that $\sum a_i
C_i$ is big and nef. The \vb $\ke_0$ is called simplifying if for two
framed bundles $\ke$ and
$\ke'$ the group $H^0(X,\kh om(\ke,\ke')(-C))$ vanishes. At the first
glance it is surprising
that there are no further stability conditions for such pairs. However,
in
many situations the general stability conditions of chapter one are hidden
behind the concept
of framed bundles.\\
\begin{definition} For $0<s<r$ the number $\nu_s(\ke_0,C_i)$ is defined as
the maximum of
${\deg(\kf)}/{s}-{\deg({\ke_0}|_{C_i})}/{r}$, where
$\kf\subset{\ke_0}|_{C_i}$ is a \vb of
rank $s$.\end{definition}
In the following we assume, that there are nonnegative integers
$a_i$, s.t. $H=\sum a_i C_i$ is
ample. This is equivalent to saying that $X\setminus C$ is affine.

\begin{proposition} If $\delta_1$ is positive with
$$\max_{0<s<r}\{{r\cdot s}/({r-s})\sum a_i\nu_s(\ke_0
,C_i)\}<\delta_1<(r-1)(C.H)~,$$ then every \vb
$\ke$ of rank $r$ together with an isomorphism
$\alpha:\ke_C\cong\ke_0$ forms a $\mu-$stable
pair $(\ke,\alpha)$.\end{proposition}
\prf The $\mu-$stability for such pairs is defined by
the following two inequalities:\\
{\it i)} ${\deg\kg}/{\rk\kg}<{\deg\ke}/{r}-{\delta_1}/{r}$
for every
\vb $\kg\subset\kea$ with $0<\rk\kg<r$ and\\
{\it ii)} ${\deg\kg}/{\rk\kg}<{\deg\ke}/{r}+\delta_1({r-\rk\kg})/(
{r\cdot\rk\kg})$
for every \vb $\kg\subset\ke$ with $0<\rk\kg<r$.\\
We first check {\it ii)}. It is enough to consider \vbs $\kg$, s.t. the
quotient $\ke/\kg$ is
torsionfree. In particular we can assume, that $\kg_{C_i}\to\ke_{C_i}$ is
injective. Then
we conclude ${\deg\kg}/{\rk\kg}={c_1(\kg).H}/{\rk\kg}=
\sum a_i{\deg(\kg_{C_i})}/{\rk
\kg}\leq\sum a_i({\deg(\ke_0)_{C_i}}/{r}+\nu_{\rk\kg}(\ke_0,C_i))=
{\deg\ke}/{r}+\sum a_i\nu
_{\rk\kg}(\ke_0,C_i)<{\deg\ke}/{r}+(({r-\rk\kg})/{r\cdot\rk\kg})\delta_1$.
To prove {\it i)}
one uses $\kea=\ke(-C)$ and {\it ii)}: For $\kg\subset\kea$ the inequality
{\it ii)}
applied to $\kg(C)\subset\ke$ implies ${\deg\kg}/{\rk\kg}+{C.H}={\deg\kg(C)}/
{\rk\kg}<{\deg\ke}/{r}+(({r-\rk\kg})/{r\cdot\rk\kg})\delta_1$.
Therefore $\delta_1<(r-1)C.H$
suffices to give {\it i)}.\qed
\begin{corollary}\label{framedareproj}
For $\max_{0<s<r}\{{r\cdot s}/({r-s})\sum a_i\nu_s(\ke_0
,C_i)\}<(r-1)(C.H)$ and $C$, such that there exists an effective,
ample divisor $H$, whose support
is contained in $C$, the moduli spaces $\km^{fr}_{X/C/\ke_0/\chi}$
of
framed \vbs are quasi-projective.\end{corollary}
\prf These \mss are in fact open subsets of the $\mu-$stable part of
the \ms
of all semistable pairs $(\ke,\alpha)$.\qed

There is a special interest in the case $\ke_0\cong\ko_C^{\oplus r}$, since
the corresponding
\mss are in fact invariants of the affine surface $X\setminus C$ (\cite{L2}).
In this case all the
numbers $\nu_s(\ke_0,C_i)$ vanish. Therefore a trivially framed bundle gives
a $\mu-$stable pair
$(\ke,\alpha)$ with respect to every $\delta_1<(r-1)C.H$.\\
In (\cite{l4},2.1.5.) a sufficient condition for a bundle $\ke_0$ to be
simplifying is
proven:
If ${\rm Hom}(\ke_0,\ke_0(-kC))=0$ for all $k>0$, then $\ke_0$ is simplifying.
We remark that at least in the rank two case this condition is closely
related to the numerical
condition we gave. It is possible to make the  condition finer, because
in the definition of the
numbers $\nu_{C_i}$ it is sufficient to take the maximum over those bundles,
which actually
live on $X$.


\end{document}